\documentclass[10pt,twoside]{article}
\usepackage[latin1]{inputenc}
\usepackage{amsmath}
\usepackage{ epsfig}
\usepackage{graphicx}
\usepackage{yhmath}
\usepackage[table]{xcolor}
\usepackage{ mathrsfs} 

\input epsf
\newcommand{\be}{\begin{equation}}
\newcommand{\ee}{\end{equation}}
\newcommand{\bd}{\begin{displaymath}}
\newcommand{\ed}{\end{displaymath}}

\newcommand{\ba}{\begin{eqnarray}}
\newcommand{\ea}{\end{eqnarray}}
\newcommand{\ban}{\begin{eqnarray*}}
\newcommand{\ean}{\end{eqnarray*}}

\newcommand{\R} {I\!\!R}
\newcommand{\E} {E\,}
\newcommand{\N} {I\!\! N}
\newcommand{\Z} {{\bf Z}}

\newcommand{\var} {\mbox{Var}\hspace{0.4mm}}

\newcommand\esp[1]{\,\,\, \text{#1}\,\,\,}
\newcommand\lct{\widehat{\lambda}_t(\theta)}
\newcommand\ftt{f_{\theta}^t}
\newcommand\fttc{\widehat{f}_{\theta}^t}
\newcommand\ltt{\ell_t(\theta)}
\newcommand\lttc{\widehat{\ell}_t(\theta)}
\newcommand\ajo{\alpha_j^{(0)}}
\newcommand\nf{E}
\newcommand\nc{\underline{c}}
\newcommand\norme[1]{\left\Vert #1 \right\Vert}
\newcommand\abs[1]{\left\vert #1 \right\vert}
\DeclareMathOperator\argmax{argmax}

\newcommand\tnc{\widehat{\theta}_n}
\newcommand\ddt{\frac{\partial}{\partial\theta}}
\newcommand\gr{>}
\newcommand\pet{<}

\newcommand\ddi{\frac{\partial}{\partial\theta_i}}
\newcommand\ddj{\frac{\partial}{\partial\theta_j}}
\newcommand\dcij{\frac{\partial^2}{\partial\theta_i\partial\theta_j}}
\newcommand\dcji{\frac{\partial^2}{\partial\theta_j\partial\theta_i}}
\renewcommand{\Box}{\hfill\rule{0.25cm}{0.25cm}} 

\newtheorem{Prop}{Proposition}[section]
\newtheorem{lem}{Lemma}[section]
\newtheorem{Theo}{Theorem}[section]

\newtheorem{rem}{Remark}[section]
\newenvironment{dem}{\ \\ {\bf Proof. }}
{\Box\par\medskip\noindent}

\def\1{{\bf 1}}

\begin{document}
\title{\bf Inference and testing for structural change in time series of counts model}
 \maketitle \vspace{-1.0cm}
\begin{center}
  Paul DOUKHAN$^*$ and William~KENGNE
\end{center}
\begin{center}
 {\it AGM, Université de Cergy Pontoise, 2 avenue Adolphe Chauvin. 95302 Cergy-Pontoise, France  (\footnote{Supported by Laboratory of Excellence MME-DII {\tt http://labex-mme-dii.u-cergy.fr/}}).\\
$^*$ IUF, Universitary institute of France.\\
  E-mail: paul.doukhan@u-cergy.fr ; william.kengne@u-cergy.fr}
\end{center}

 \pagestyle{myheadings}
 \markboth{Inference and testing for structural change in time series of counts model}{P. Doukhan and W. Kengne}

\textbf{Abstract} : We  consider here together the inference questions and the change-point problem
in Poisson autoregressions (see Tj{\o}stheim, 2012).
 The conditional mean (or intensity)  of the process is involved as a non-linear function of it past values and the past observations.
 Under Lipschitz-type conditions, it is shown that the conditional mean can be written as a function of lagged observations.
 In the latter model,  assume that the link function depends on an unknown parameter $\theta_0$. The consistency and the asymptotic normality
 of the maximum likelihood estimator of the parameter are proved. These results are used to study change-point problem in the parameter $\theta_0$.
 We propose two tests  based on the likelihood of the observations. Under the null hypothesis (i.e. no change), it is proved that both  those test
  statistics converge to an explicit distribution. Consistencies under  alternatives are proved for both tests. Simulation results  show how those procedure work  practically, and an application to real data is also processed. ~\\ ~\\


 \section{Introduction} \label{intro}
  Count models are of a large current interest, see the discussion paper  Tj{\o}stheim (2012).
 Integer-valued time series appear as natural for modeling count events. Examples may be found in epidemiology (number of new infections),
  in finance (number of transactions per minute), in industrial quality control (number of defects);
  see for instance Held {\it et al.} (2005),  Brännäs and Quoreshi (2010) and  Lambert (1992).
\\
Real data are definitely not stationary. Several ways to consider such structural changes are possible as this was demonstrated during  the thematic cycle {\it Nonstationarity and Risk Management} held in Cergy-Pontoise during year 2012 (\footnote{\tt http://www.u-cergy.fr/en/advanced-studies-institute/thematic-cycles/thematic-cycle-2012/finance-cycle.html}).
Structural breaks are a reasonable possibility for this, when no additional knowledge is available. The paper is thus aimed at considering such changes in regime for a large class of integer valued time series.\\
  Let $Y=(Y_t)_{t \in \Z}$ be an integer-valued time series and
   $\mathcal{F}_t=\sigma(Y_s, ~ s \leq t)$ the $\sigma$-field generated by the whole past at time $t$, we denote by $\mathcal{L}(Y_t/\mathcal{F}_{t-1})$ the
 conditional  distribution of $Y_t$ given the past. Models with various marginal distributions and dependence
   structures have been studied for instance Kedem and Fokianos (2002), Davis {\it et al.} (2005), Ferland {\it et al.} (2006), Davis and Wu (2009),
   Wei$ß$ (2009).
\\
   Fokianos {\it et al.} (2009) considered the Poisson autoregression such that $\mathcal{L}(Y_t/\mathcal{F}_{t-1})$ is Poisson distributed with a intensity $\lambda_t$
   which is a function of $\lambda_{t-1}$ and $Y_{t-1}$.  Under linear autoregression, they prove both the consistency and the asymptotic normality
   of the maximum likelihood estimator of the regression parameter, by using a perturbation approach, which allows to used the standard Markov ergodic setting.
    Fokianos and Tjøstheim (2012) extend the method
   to nonlinear Poisson autoregression with $\lambda_t=f(\lambda_{t-1}) + b(Y_{t-1})$ for nonlinear  measurable functions $f$ and $b$.
   In the same vein, Neumann (2011) had studied some large situation where $\lambda_t=f(\lambda_{t-1},Y_{t-1})$. He focused on the absolute regularity and the
   specification test for the intensity function, while the recent work of Fokianos and  Neumann (2013) studied goodness-of-fit tests which
   are able to detect local alternatives.
    Doukhan {\it et al.} (2012), consider  a more general model with infinitely many lags.
   Stationarity  and the existence of moments are proved   by using a weak dependent approach and contraction properties.\\
  Later, Davis and Liu (2012) studied the model where the distribution $\mathcal{L}(Y_t/\mathcal{F}_{t-1})$ belongs to a class of
   one-parameter exponential family with finite order dependence. This class contains Poisson and negative binomial distribution.
   From the theory of iterated random functions, they establish the stationarity and the absolute regularity properties of the process.
   They also prove the consistency and  asymptotic normality of the maximum likelihood estimator of the parameter of the model.
   \\
    Douc {\it et al.} (2012) considers a class of observation-driven time series which covers linear, log-linear, and threshold Poisson autoregressions.
   They approach is based on a recent  theory for Markov chains based upon Hairer and Mattingly (2006) recent work; this allows existence and uniqueness of the invariant distribution for Markov chains without irreducibility. They also proved  consistency of the conditional likelihood estimator of the model (even under false models); the asymptotic normality is not yet considered in this setting.
   \\
   For change point procedures,   Kang and Lee (2009) proposed CUSUM procedure for testing for parameter change in a first-order random coefficient integer-valued autoregressive model
      defined through thinning operator. Fokianos and Fried (2010, 2012) studied mean shift in linear and log-linear Poisson autoregression. The dependence between the level shift and time
      allows their model to detect several types of interventions effects such as outliers. Franke  {\it et al.} (2012) consider parameter change in Poisson
      autoregression of order one. Their test are based on the cumulative sum of residuals using conditional least-squares estimator.\\
\medskip

\noindent
In this article, we first consider a time series of counts $Y=(Y_t)_{t \in \Z}$ satisfying :
   \begin{equation}\label{model1}
   Y_t/\mathcal{F}_{t-1}   \sim\mbox{ Poisson}(\lambda_t) ~ \text{with} ~ \lambda_t = F(\lambda_{t-1}, \ldots; Y_{t-1}, \ldots)
   \end{equation}
   where  $\mathcal{F}_t=\sigma(Y_s, ~ s \leq t)$  and $F$ a measurable non-negative function.
   The properties of the general class of Poisson autoregression model (\ref{model1}) have been studied by Doukhan {\it et al.} \cite{Doukhan2012}.
   According to the fact that the model (\ref{model1}) can take into account the whole past information of the process, it dependence structure is more
   general than those studied before. Proceeding as in Doukhan and Wintenberger (2008), we show that under some Lipschitz-type conditions on $F$,
    the conditional mean $\lambda_t$ can be written as a function of the past observations. This leads us to consider the model
   \begin{equation}\label{model2}
   Y_t/\mathcal{F}_{t-1}  \sim \mbox{ Poisson}(\lambda_t) ~ \text{with} ~ \lambda_t = f(Y_{t-1}, \ldots)
   \end{equation}
   where $f$ is a measurable non-negative function. We assume that $f$ is know up to some parameter $\theta_0$ belonging to a compact set $\Theta$.
    That is
     \begin{equation}\label{model3}
   Y_t/\mathcal{F}_{t-1}  \sim \mbox{ Poisson}(\lambda_t) ~ \text{with} ~ \lambda_t = f_{\theta_0}( Y_{t-1}, \ldots) ~ \esp{and} ~ \theta_0 \in \Theta.
   \end{equation}
    Many classical integer-valued time series satisfying (\ref{model3}) (see the examples below).\\
    Remark that, model (\ref{model3}) (as well as models (\ref{model1}) and (\ref{model2})) can be represented in terms of Poisson processes.
    Let $(N_t)_{t\geq 0}$ be a sequence of independent Poisson processes of unit intensity. $Y_t$ can be seen as a number (say $N_t(\lambda_t)$) of events
    that occurred in the time interval $[0, \lambda_t ]$. So, we have the representation
     \begin{equation}\label{model3bis}
   Y_t = N_t(\lambda_t) ~ \text{with} ~ \lambda_t = f_{\theta_0}( Y_{t-1}, \ldots).
   \end{equation}

     The paper is first work out the asymptotic properties of the maximum likelihood estimator of the model (\ref{model3}).
     We assume that the function $f_{\theta}$ satisfies some classical Lipschitz-type inequality and investigate sufficient conditions
     for the consistency and the asymptotic normality of the maximum likelihood estimator of $\theta_0$. Contrary to Fokianos {\it et al.} \cite{Fokianos2012}
     the increasing of the function $(y_k)_{k\in \N} \mapsto f_{\theta}((y_k)_{k\in \N})$  (not easy to define here) as well as the existence of the fourth order
     derivative of the function $\theta \mapsto f_{\theta}$ are not needed. Although the models  studied by Davis and Liu \cite{Davis2012} and
     Douc {\it et al.} \cite{Douc2012} allow large classes of marginal distributions, the infinitely  many lags of model (\ref{model1}) (or model (\ref{model3}))
     makes it allows a large type of dependence structure.

     The second contribution of this work are the two tests for change detection in model (\ref{model3}). We propose a new idea to take into account the change-point
     alternative. This make that, the procedures proposed will be numerically easy to apply than that proposed by Kengne (2012). The consistency under the alternative
     are proved. Contrary to Franke  {\it et al.} \cite{Franke2012}, the multiple change alternative has been considered and the independence between
     the observation before and after the change-point is not assumed. Note that, the intervention problem studied by Fokianos and Fried
     \cite{Fokianos2010, Fokianos2012_itv} is intended to sudden shift in the conditional means of the process. In the classical change-point setting, such
      intervention will be asymptotically negligible.\\

      In the forthcoming Section \ref{asum}, it is provided some assumptions on model with examples. Section \ref{infer} is devoted to the definition of the maximum likelihood estimator
      with it asymptotic properties. In Section \ref{change}, we propose the tests for detecting change in parameter of model (\ref{model3}). Some simulation results
      and real data application are presented in Sections \ref{simu} and \ref{rdata} and the proofs of the main results are provided in Section \ref{proof}.

    \section{Assumptions and examples}\label{asum}

    \subsection{Assumptions}
   We will use the following classical notations:
 \begin{enumerate}
    \item $\|y\|:=\sum\limits_{j=1}^{p} |y_j|$ for any $y \in \R^p$;
    \item for any compact set $ \mathcal{K} \subseteq\R^d$ and for any function
    $g:\mathcal{K} \longrightarrow\R^{d'}$, $ \|g\|_\mathcal{K} =\sup_{\theta\in \mathcal{K} }(\|g(\theta)\|)$;
    \item it $Y$ is a random variable with finite moment of order $r>0$, then $\|Y\|_r = (\E|Y|^r)^{1/r} $;
    \item for any set $ \mathcal{K} \subseteq\R^d$,  $\overset{\circ}{\mathcal{K} }$ denotes the interior of $\mathcal{K} $.
 \end{enumerate}

 ~ \\
    \noindent A classical Lipschitz-type conditions is assumed on the model (\ref{model1}). \\

 \noindent {\bf Assumption A$_F$} : There exists a sequence of non-negative real numbers $(\alpha_j)_{j\geq 1}$ satisfying
   $ \sum\limits_{j=1}^{\infty} \alpha_j <1$ and such that for any $y,y' \in ((0,\infty)\times\N)^{\N}$,
\begin{equation*}
 |F(y)-F(y')| \leq \sum\limits_{j=1}^{\infty}\alpha_j \|y_j-y'_j\|.
\end{equation*}
  Under assumption (\textbf{A}$_F$), Doukhan {\it et al.} \cite{Doukhan2012} have shown that there exists a solution $(Y_t,\lambda_t)_{t\in \Z}$ which is
  $\tau$-weakly dependent strictly stationary with finite moment of any order.
  The following proposition show that the conditional mean $\lambda_t$ of model (\ref{model1}) can be expressed as a function of only the past observations
  of the process.
 \begin{Prop}\label{prop1}
 Under (\textbf{A}$_F$), there exists a measurable function $f: \R^+ \rightarrow \R^+ $ such that the strictly stationary ergodic solution of model (\ref{model1})
 satisfying almost surely
 $$ \lambda_t = f(Y_{t-1},\dots).$$
 \end{Prop}
 \noindent Proposition \ref{prop1} shows that, the information contained in the unobservable process $(\lambda_{t-j})$ can be capture by the observable process $(Y_{t-j})$.
 This shows that, in practice the autoregression on $(Y_{t-j})$ is sufficient to capture the information of the whole past.
 This representation is very important in the inference framework. Note that, if one carries inference on the model (\ref{model1}) by assuming that
 $\lambda_t = F_{\theta_0}(\lambda_{t-1}, \ldots; Y_{t-1}, \ldots) $, it will not be easy to compute $\partial \lambda_t/\partial \theta $ or to write it as a function
 of $\partial F_{\theta} /\partial \theta $. So, the asymptotic normality of the maximum likelihood estimator of $\theta_0$ (see (\ref{loglik}) and (\ref{emv}))
 will be very difficult to study. \\
 \indent We focus on the model (\ref{model3}) with the following assumptions. For $i=0,\, 1, \, 2$ and for any compact set $\mathcal{K} \subseteq \Theta$,
 define \\
 {\bf Assumption A$_i(\mathcal{K} )$}: $\|{\partial^i f_\theta(0)}/{\partial\theta^i}\|_\Theta<\infty$
 and there exists a sequence of non-negative real numbers $(\alpha^{(i)}_k(\mathcal{K} ))_{k\geq 1}$ satisfying
   $ \sum\limits_{j=1}^{\infty} \alpha^{(0)}_k(\mathcal{K} ) < 1$ or  $ \sum\limits_{j=1}^{\infty} \alpha^{(i)}_k(\mathcal{K} ) <\infty$ (for $i=1,2$) such that
\begin{equation*}
 \Big\|\dfrac{\partial^i f_\theta(y)}{\partial\theta^i}-\dfrac{\partial^i f_\theta(y')}{\partial\theta^i}\Big\|_\mathcal{K}
 \leq \sum\limits_{k=1}^{\infty}\alpha^{(i)}_k(\mathcal{K})|y_k-y'_k|\quad \mbox{ for all } y, y' \in (\R^+)^{\N}.
\end{equation*}

 \noindent The Lipschitz-type condition A$_0(\Theta)$ is the version of model (\ref{model3}) of the assumption  \textbf{A}$_F$.
  It is classical when studying the existence of solutions of such model (see for instance \cite{Doukhan2008}, \cite{Bardet2009} or \cite{Doukhan2012}).
  The assumptions A$_1(\mathcal{K})$ and A$_2(\mathcal{K})$ as well as the following assumptions D$(\Theta)$, Id($\Theta$) and  Var($\Theta$)
  are needed to define and to study the asymptotic properties of the maximum likelihood estimator of the model (\ref{model3}).\\

  \noindent {\bf Assumption D$(\Theta)$:} $\exists\underline{c}>0$ such that
$\displaystyle \inf_{\theta \in
 \Theta}(f_\theta(y))\geq \underline{c}$ for all $y\in (\R^+)^{\N}.$ \\
~\\
 {\bf Assumption Id($\Theta$):} For all  $(\theta,\theta')\in \Theta^2$,
 $ \Big( f_{\theta}(Y_{t-1},\dots)=f_{\theta'}(Y_{t-1},\dots)  \ \text{a.s.} ~ \text{ for some } t \in \Z \Big) \Rightarrow \ \theta = \theta'.$
 {\bf Assumption Var($\Theta$):} For all $\theta  \in \Theta $ and $t \in \Z$,
 the components of the vector $\dfrac{\partial f_{\theta}}{\partial \theta^{i}}(Y_{t-1,\dots})$  are  a.s. linearly independent.

  \subsection{Examples}

  \subsubsection{Linear Poisson autoregression}

   We consider an integer-valued time series $(Y_t)_{t \in \Z}$ satisfying for any $t \in \Z$
   \begin{equation}\label{ex1}
   Y_t/\mathcal{F}_{t-1}  \sim \mbox{ Poisson}(\lambda_t) ~ \text{with} ~ \lambda_t = \phi_0(\theta_0) + \sum_{k \geq 1} \phi_k(\theta_0) Y_{t-k}
   \end{equation}
 where $\theta_0 \in \Theta \subset \R ^d$, the function $\theta \mapsto \phi_k(\theta)$ are positive and satisfying $\sum_{k \geq 1} \norme{\phi_k(\theta)}_\Theta < 1$.
 Assumptions A$_0(\Theta)$ holds automatically. If the function $\phi_k$ are twice continuous differentiable such that
 $\sum_{k \geq 1} \norme{\phi_k '(\theta)}_\Theta < \infty$ and $\sum_{k \geq 1} \norme{\phi_k '' (\theta)}_\Theta < \infty$ then A$_1(\Theta)$ and A$_2(\Theta)$
 hold. If $ \underset{\theta \in \Theta} {\inf} \phi_0(\theta) >0$ then D$(\Theta)$ holds.
 Moreover, if there exists a finite subset $I \subset \N - \{0\}$ such that the function $\theta \mapsto (\phi_k(\theta))_{ k \in I}$ is injective, then assumption
 Id$(\Theta)$ holds i.e. model (\ref{ex1}) is identifiable. \\
 Note that the popular Poisson INGARCH model (see \cite{Ferland2006} or \cite{Weib2009})  is a special case of model (\ref{ex1}). Finally, the model (\ref{ex1})
 can be generalized by considering
      \[\lambda_t = \phi_0(\theta_0) + \sum_{k \geq 1} \phi_k(\theta_0) h(Y_{t-k}) \]
      where $h$ is a Lipschitz function. 
 The threshold model can be obtained with $h(y)= (y-\ell)_+$ for some $\ell >0$, where $y_+ = \max(y,0)$.

  \subsubsection{Power Poisson autoregression}
  We consider a power Poisson INGARCH($p,q$) process defined by
   \begin{equation}\label{ex2}
   Y_t/\mathcal{F}_{t-1}  \sim  \mbox{ Poisson}(\lambda_t) ~ \text{with} ~ \lambda_t = \Big( \alpha_0(\theta_0) + \sum_{k = 1}^{p} \alpha_k(\theta_0)\lambda_{t-k}^{\delta}
    + \sum_{k = 1}^{q} \beta_k(\theta_0)Y_{t-k}^{\delta} \Big)^{1/\delta}
   \end{equation}
  where $\delta \geq 1$, $\theta_0 \in \Theta \subset \R ^d$ and the function $\theta \mapsto \alpha_k(\theta)$ and $\theta \mapsto \beta_k(\theta)$ are positive.
  If
  $ \sum_{k = 1}^{p} \norme{\alpha_k(\theta)}_\Theta^{1/\delta} + \sum_{k = 1}^{q} \norme{\beta_k(\theta)}_{\Theta}^{1/\delta} <1 $, then the Lipschitz-type
  condition (\textbf{A}$_F$) is satisfied. In this case, we can find a sequence of non-negative function $(\psi_k(\theta))_{k \geq 1}$ such that
 $\lambda_t = \Big( \psi_0(\theta_0) + \sum_{k \geq 1} \psi_k(\theta_0)Y_{t-k}^{\delta} \Big)^{1/\delta}$.
 If $ \underset{\theta \in \Theta} {\inf} \alpha_0(\theta) >0$ then D$(\Theta)$ holds.
 Moreover, if there exists a finite subset $I \in \N - \{0\}$ such that the function $\theta \mapsto (\psi_k(\theta))_{ k \in I}$ is injective, then assumption
 Id$(\Theta)$ holds i.e. model (\ref{ex2}) is identifiable. \\

 \section{Likelihood inference}\label{infer}

Assume that the trajectory $(Y_1, \dots, Y_n)$ is observed. The conditional (log)-likelihood
(up to a constant) of model~$(3)$ computed on $T \subset\{1,\dots,n\}$, is given by
 \[ L_n(T,\theta) = \sum_{t \in T}(Y_t\log \lambda_t(\theta)- \lambda_t(\theta)) = \sum_{t\in T} \ell_t(\theta) \esp{with}
  \ell_t(\theta) = Y_t\log \lambda_t(\theta)- \lambda_t(\theta)\]
  where $ \lambda_t(\theta) = f_{\theta}(Y_{t-1},\dots)$. In the sequel, we use the notation $f^t_{\theta} := f_{\theta}(Y_{t-1}, \ldots).$
 Since only $Y_1,\dots,Y_n$ are observed, we will use an approximation version of the (log)-likelihood defined by
 \begin{equation}\label{loglik}
 \widehat{L}_n(T,\theta) = \sum_{t \in T}(Y_t\log\widehat{\lambda}_t(\theta)-\widehat{\lambda}_t(\theta)) = \sum_{t\in T}\widehat{\ell}_t(\theta) \esp{with}
 \widehat{\ell}_t(\theta) = Y_t\log\widehat{\lambda}_t(\theta)-\widehat{\lambda}_t(\theta)
 \end{equation}
 where
 $\lct := \widehat{f}^t_{\theta}:= f_{\theta}(Y_{t-1}, \dots, Y_1,0,\dots)$ and $\widehat{\lambda}_1(\theta) = f_{\theta}(0,\dots) $.
 The maximum likelihood  estimator (MLE) of $\theta_0$ computed on $T$ is defined by
 \begin{equation}\label{emv}
  \widehat{\theta}_n(T) = \argmax_{\theta \in \Theta} (\widehat{L}_n(T,\theta)).
  \end{equation}
 For any $k,k' \in \Z$ such as $k\leq k'$, denote
 \[T_{k,k'} = \{k,k+1, \ldots,k' \} .\]

 \begin{Theo}\label{theo1}
 Let $(j_n)_{n\geq 1}$ and $(k_n)_{n\geq 1}$ be two integer valued sequences such that
 $j_n \leq k_n$, $k_n \to +\infty$ and $k_n-j_n \to +\infty$ as $n\to +\infty$. Assume  $\theta_0 \in  \ring{\Theta}$ and $D(\Theta)$,
 $\mathrm{Id}(\Theta)$ and $A_0(\Theta)$ hold with
 \begin{equation} \label{theo1_eq}
 \sum_{j \geq 1}\log j\times \alpha_j^{(0)}(\Theta) < \infty.
 \end{equation}
 It holds that
 \[\widehat{\theta}_n(T_{j_n,k_n}) \xrightarrow [n\to +\infty]{a.s.} \theta_0.\]

 \end{Theo}

 \noindent The following theorem shows the asymptotic normality of the MLE of model (\ref{model3}).

 \begin{Theo}\label{theo2}
 Let $(j_n)_{n\geq 1}$ and $(k_n)_{n\geq 1}$ be two integer valued sequences such that
 $j_n \leq k_n$, $k_n \to +\infty$ and $k_n-j_n \to +\infty$ as $n\to +\infty$. Assume  $\theta_0 \in  \ring{\Theta}$ and $D(\Theta)$,
 $\mathrm{Id}(\Theta)$,  Var($\Theta$) and $A_i(\Theta)$ $i=0,1,2$ hold with
 \begin{equation} \label{theo2_eq}
 \sum_{j \geq 1} \sqrt{j} \times \alpha_j^{(i)}(\Theta) < \infty.
  \end{equation}
 It holds that
 \[ \sqrt{k_n-j_n}(\widehat{\theta}_n(T_{j_n,k_n}) - \theta_0) \xrightarrow [n\to +\infty]{\mathcal{D}} \mathcal{N}(0,\Sigma^{-1} )\]
 where $\Sigma=E\big(\frac{1}{f_{\theta_0}^0}(\ddt f^0_{\theta_0})(\ddt f^0_{\theta_0})'\big).$

 \end{Theo}
 ~ \\

 \noindent According to the Lemma \ref{lem1} and the proof of Theorem \ref{theo2}, the matrix
 \[ \Big(\frac{1}{k_n-j_n} \sum_{t=j_n}^{k_n}\frac{1}{\fttc}\big(\frac{\partial}{\partial\theta}\fttc\big)\big(\frac{\partial}{\partial\theta}\fttc\big)'\Big)\Big|_{\theta = \widehat{\theta}_n(T_{j_n,k_n})}
  \esp{and} \Big(- \dfrac{1}{k_n-j_n} \dfrac{\partial^2}{\partial \theta \partial \theta ' } \widehat{L}_n(T_{j_n,k_n},\theta) \Big)\Big|_{\theta = \widehat{\theta}_n(T_{j_n,k_n})} \]
  are consistent estimators of the covariance matrix $\Sigma$.

\begin{rem}
    \begin{enumerate}
     \item In  Theorems \ref{theo1} and  \ref{theo2}, the typical sequences $j_n=1$ and $k_n=n$, $\forall n \geq 1$ can be chosen.
           This choice is the case where the estimator is computed with all the observations. But in the change-point study, the estimator needs to be
           calculated over a part of the observations. Results are written this way to cover the change point situation.
     \item If the Lipschitz coefficients  $(\alpha_j^{(i)}(\Theta))_{j\geq 1}$ satisfy $\alpha_j^{(i)}(\Theta) = \mathcal{O}(j^{-\gamma})$ with $\gamma >1/2$, then
      the conditions (\ref{theo1_eq}) and  (\ref{theo2_eq}) of Theorem \ref{theo1} and Theorem \ref{theo2} hold.

    \end{enumerate}
\end{rem}

 \section{Testing for Parameter Changes}\label{change}
  We consider the observations $Y_1, \dots, Y_n$ generated as in model $(3)$ and assume that the parameter  $\theta_0$ may
 change over time. More precisely, we assume that
 $\exists K \geq 1, \theta_1^*, \dots, \theta_K^* \in \Theta, 0 = t_0^* \pet t_1^* \pet \dots \pet t_{K-1}^* \pet t_K^* = n$ such that $Y_t = Y_t^{(j)}$ for
 $t^*_{j-1} \pet t \leq t^*_j$, where the process $(Y_j^{(j)})_{t\in \Z}$ is a stationary solution of $(3)$ depending on $\theta^*_j$.
 The case where the parameter does not change corresponds to $K = 1$. For this problem, we consider the following hypothesis:
\begin{enumerate}
  \item[H$_0$:] The observations $(Y_1, \dots, Y_n)$ are a trajectory of a process $(Y_t)_{t\in \Z}$ solution of~$(3)$, depending on
  $\theta_0 \in \Theta$.
  \item[H$_1$:] There exists  $K \geq 2, \theta_1^*, \theta_2^*, \dots, \theta_K^*$ with $\theta_1^*\neq \theta_2^*\neq \dots \neq \theta_K^*,
  0 = t_0^* \pet t_1^* \pet \cdots \pet t_{K-1}^* \pet t_K^* = n$ such that the observations $(Y_t)_{t^*_{j-1} \pet t \leq t^*_j}$
  are a trajectory of the process $(Y_t^{(j)})_{t \in \Z}$ solution of~$(3)$, depending on $\theta^*_j$.
\end{enumerate}
 Let us note that, contrary to Franke {\it et al.} (2012), the independence between the observations before and after the change-point is not assumed.
 Moreover, we can consider more complex parameter in the model, not only the constant term of the conditional mean as in \cite{Franke2012}. In the case of linear
 Poisson autoregression, the assumption (A9) of their model leads to the change in the unconditional mean. Therefore, our change-point problem is more general.
 \begin{rem}
   Let $\lambda^{(j)}_t$ be the conditional mean on the segment $T^*_j=\{t^*_{j-1},t^*_{j-1}+1, \ldots,t^*_{j}\}$.
   By using the representation (\ref{model3bis}), we can write
  \[  Y_t^{(j)} = N_t(\lambda_t^{(j)}) \esp{with} \lambda^{(j)}_t = f_{\theta^*_j}(Y_{t-1}^{(j)}, \ldots).\]
   We can also write
   \[  Y_t^{(j)} =\tilde{F}(Y_{t-1}^{(j)}, \ldots; N_t)  \esp{where} \tilde{F}(y_1,y_2, \ldots; N_t) = N_t(f_{\theta^*_j}(y_1,y_2, \ldots))
   \esp{for any} y_k \in \N, ~ k\geq 1 .\]
   So, for any $y=(y_k)_{k\geq1}$ and $y'=(y'_k)_{k\geq1}$, we have
    \begin{equation} \label{approx_F_tilde}
    \E |\tilde{F}(y; N_t) - \tilde{F}(y'; N_t)| = \E |N_t(f_{\theta^*_j}(y)) - N_t(f_{\theta^*_j}(y'))| = |f_{\theta^*_j}(y) - f_{\theta^*_j}(y')|
       \leq \sum_{k\geq 1} \alpha^{(0)}_k |y_k-y'_k|,
    \end{equation}
  where the second equality follows by seen $|N_t(f_{\theta^*_j}(y)) - N_t(f_{\theta^*_j}(y'))|$ as a number of events $N_t$ that occur in the time interval
  $[0, |f_{\theta^*_j}(y) - f_{\theta^*_j}(y')|]$.\\
  Let $(\widetilde{Y}_t,\widetilde{\lambda}_t)_{t^*_{j-1} \pet t \leq t^*_j}$ be the nonstationary approximation of the process
  $(Y_t^{(j)},\lambda_t^{(j)})_{t \in \Z}$ on the segment $T^*_j$; i.e.
   \[  \widetilde{Y}_t = N_t(\widetilde{\lambda}_t) \esp{with}
   \widetilde{\lambda}_t = f_{\theta^*_j}(Y_{t-1}^{(j)}, \ldots, Y_{t-t^*_{j-1}}^{(j)},Y_{t-t^*_{j-1}-1}^{(j-1)}, \ldots,Y^{(1)}_1,0, \ldots).\]
  By using assumption A$_0(\Theta)$ and relation (\ref{approx_F_tilde}), one can show that the approximated process
  $(\widetilde{Y}_t,\widetilde{\lambda}_t)_{t^*_{j-1} \pet t \leq t^*_j}$ converges (in $L^r$ for any $r\geq1$) to the stationary regime
  (see for instance Bardet {\it et al.} \cite{Bardet2010} where similar apprimation has been done in the case of causal time series).
  So, the results of the Section \ref{asymp_H1} may be extended (modulo the validity of the approximation) by relaxing the stationarity assumption after change.

 \end{rem}

 Recall that under H$_0$, the likelihood function of the model computed on  $T \subset \{1, \dots, n\}$ is given by
 \[\widehat{L}_n(T,\theta) = \sum_{t\in T}(Y_t\log\fttc-\fttc)\]
 where $\fttc = f_{\theta}(Y_{t,1}, \dots)$ and the maximum likelihood estimator is given by $\widehat{\theta}_n (T) = \argmax \widehat{L}_n(T, \theta)$.
 It holds from Theorem \ref{theo2} that, under H$_0$, the asymptotic covariance matrix of
 $\widehat{\theta}_n(T_{1,n})$ is $\widehat{\Sigma}_n^{-1}$ where
 \[\widehat{\Sigma}_n = \Big(\frac{1}{n} \sum_{t=1}^n\frac{1}{\fttc}\big(\frac{\partial}{\partial\theta}\fttc\big)\big(\frac{\partial}{\partial\theta}\fttc\big)'\Big)\Big|_{\theta = \widehat{\theta}_n(T_{1,n})}.\]
 $\widehat{\Sigma}_n$ is a consistant estimator of
\[\Sigma = E\Big(\frac{1}{f^0_{\theta_0}}\big(\frac{\partial}{\partial\theta}f^0_{\theta_0}\big)\big(\frac{\partial}{\partial\theta}f^0_{\theta_0}\big)'\Big)\]
 under H$_0$ (see the proof of Theorem \ref{theo2}). The consistency of $\widehat{\Sigma}_n$ under H$_1$ is not ensured.
 $\widehat{\Sigma}_n$ does not take into account the change-point alternative. So, the consistency under H$_1$ of any test based
 on $\widehat{\Sigma}_n$ will not be easy to prove.

 Let $(u_n)_{n \geq 1}$ and $(u_n)_{n \geq 1}$ be two integer valued sequences satisfying $u_n, v_n \to +\infty, \frac{u_n}{n}, \frac{v_n}{n} \to 0$ as $n \to +\infty$.
 Our test statistic  is based on the following matrix
%
  \[ \widehat{\Sigma}_n(u_n) = \frac{1}{2}\Big[\Big(\frac{1}{u_n}\sum_{t=1}^{u_n}\frac{1}{\fttc}\big(\frac{\partial}{\partial\theta}\fttc\big)
 \big(\frac{\partial}{\partial\theta}\fttc\big)'\Big)\Big|_{\theta = \widehat{\theta}_n(T_{1,u_n})} +  \Big(\frac{1}{n-u_n}\sum_{t=u_n+1}^{n}\frac{1}{\fttc}\big(\frac{\partial}{\partial\theta}\fttc\big)
 \big(\frac{\partial}{\partial\theta}\fttc\big)'\Big)\Big|_{\theta = \widehat{\theta}_n(T_{u_n+1,n})}\Big].\]
 Theorem \ref{theo1} and Lemma \ref{lem1} show that $\widehat{\Sigma}_n(u_n)$ is consistent under H$_0$. Under H$_1$, we will use
 the classical assumption that the breakpoint gown at rate~$n$. This will allow us to show that the first component
 of $\widehat{\Sigma}_n(u_n)$ converges to the covariance matrix of the stationary model of the first regime. It will
 be a key to prove the consistency under H$_1$.\\
 Another way to deal is to consider the matrix
 \[ \widetilde{\Sigma}_n(u_n) = \frac{1}{2}\Big[\Big(\frac{1}{n-u_n}\sum_{t=1}^{n-u_n}\frac{1}{\fttc}\big(\frac{\partial}{\partial\theta}\fttc\big)
 \big(\frac{\partial}{\partial\theta}\fttc\big)'\Big)\Big|_{\theta = \widehat{\theta}_n(T_{1,n-u_n})} +  \Big(\frac{1}{u_n}\sum_{t=n-u_n+1}^{n}\frac{1}{\fttc}\big(\frac{\partial}{\partial\theta}\fttc\big)
 \big(\frac{\partial}{\partial\theta}\fttc\big)'\Big)\Big|_{\theta = \widehat{\theta}_n(T_{n-u_n+1,n})}\Big].\]
 Asymptotically, both the matrices $\widehat{\Sigma}_n(u_n)$ and $\widetilde{\Sigma}_n(u_n)$ have the same behavior under H$_0$. In the case of non stationarity
 after change, the procedure using $\widetilde{\Sigma}_n(u_n)$ can provide more distortion; because due to the dependance the second component of
 $\widetilde{\Sigma}_n(u_n)$ will converge very slowly than the first component of $\widehat{\Sigma}_n(u_n)$.

 Let us define now the tests statistics:
 \begin{itemize}
  \item
 $\displaystyle \widehat{C}_n = \max_{v_n \leq k \leq n-v_n}\widehat{C}_{n,k}$ where
 \[\widehat{C}_{n,k} = \frac{1}{q^2(\frac{k}{n})}\frac{k^2(n-k)^2}{n^3}\big(\widehat{\theta}_n(T_{1,k})-\widehat{\theta}_n(T_{k+1,n})\big)'
 \widehat{\Sigma}_n(u_n)\big(\widehat{\theta}_n(T_{1,k})-\widehat{\theta}_n(T_{k+1,n})\big);\]
 where $q$ is a weight function define on $(0,1)$, see bellow;
  \item
 $\displaystyle \widehat{Q}_n^{(1)} = \max_{v_n \leq k \leq n-v_n}\widehat{Q}_{n,k}^{(1)}$ where
 \[\widehat{Q}_{n,k}^{(1)} = \frac{k^2}{n}\big(\widehat{\theta}_n(T_{1,k})-\widehat{\theta}_n(T_{1,n})\big)'
 \widehat{\Sigma}_n(u_n)\big(\widehat{\theta}_n(T_{1,k})-\widehat{\theta}_n(T_{1,n})\big);\]
  \item
 $\displaystyle \widehat{Q}_n^{(2)} = \max_{v_n \leq k \leq n-v_n}\widehat{Q}_{n,k}^{(2)}$ where
 \[\widehat{Q}_{n,k}^{(2)} = \frac{(n-k)^2}{n}\big(\widehat{\theta}_n(T_{k+1,n})-\widehat{\theta}_n(T_{1,n})\big)'
 \widehat{\Sigma}_n(u_n)\big(\widehat{\theta}_n(T_{k+1,n})-\widehat{\theta}_n(T_{1,n})\big).\]
 \end{itemize}
 The first procedure is based on the statistique $\widehat{C}_n$ and the other one is based on the statistic $\widehat{Q}_n$ defined by
 \[ \widehat{Q}_n := \max(\widehat{Q}_n^{(1)} , \widehat{Q}_n^{(2)}) .\]
 The weight function $q$ is used to increase the power of the test based on the statistic $\widehat{C}_n$. In the sequel, we will assume that\\
 $q : (0,1) \rightarrow (0,\infty) $ is  non-decreasing in a neighborhood of zero, non-increasing in a neighborhood of one and satisfying
          $\underset{\eta < \tau < 1-\eta } {\inf} q(\tau) > 0 $  for all $0< \eta < 1/2 $.\\
 The behavior of the weighted function $q$ can be controlled at the neighborhood of zero and one by the integral
 \[ I_{0,1}(q,c) = \int_0^1 \dfrac{1}{t(1-t)} \exp\Big( - \dfrac{c q^2(t)}{t(1-t)} \Big) dt , ~ c>0 \]
 see  Csörgo {\it et al.} \cite{Csorgo1986} or Csörgo and Horváth \cite{Csorgo1993}.
 The natural weight choice is $q(t) = \big( t(1-t) \big)^{\gamma}$ with $0\leq \gamma < 1/2$.\\
 Furthermore, in practice the sequences  $(u_n)_{n \geq 1}$ and $(u_n)_{n \geq 1}$ are chosen to ensure the convergence of the numerical algorithm used to compute
 $\widehat{\theta}_n(T_{1,u_n})$ and $\widehat{\theta}_n(T_{1,v_n})$. For Poisson INGARCH model, $u_n=v_n=[(\log n )^{\delta_0}]$ (with $5/2 \leq \delta_0 \leq 3$)
 can be chosen (see also  Remark 1 of  \cite{Kengne2011}).

\subsection{Asymptotic behavior under the null hypothesis}

   The asymptotic distributions  of these statistics under H$_0$ are given in the next theorem.

  \begin{Theo}\label{theo3}
 Assume  $D(\Theta)$, $\mathrm{Id}(\Theta)$,  Var($\Theta$) and $A_i(\Theta)$ $i=0,1,2$ hold with
\[ \sum_{j \geq 1} \sqrt{j} \times \alpha_j^{(i)}(\Theta) < \infty.\] Under H$_0$ with $\theta \in  \ring{\Theta}$,
\begin{enumerate}
  \item if $\exists c \gr 0$ such that $I(q,c) \pet \infty$, then
  \[\widehat{C}_n \xrightarrow[n \to +\infty]{\mathcal{D}} \sup_{0\pet \tau \pet 1}\frac{\norme{W_d(\tau)}^2}{q^2(\tau)};\]
  \item for $j = 1,2$,
  \[\widehat{Q}_n^{(j)}\xrightarrow[n\to +\infty]{\mathcal{D}}\sup_{0 \pet \tau \pet 1}\norme{W_d(\tau)}^2,\]
  where $W_d$ is a $d$-dimensional Brownian bridge.
\end{enumerate}
\end{Theo}
 The distribution of $ \underset{ 0 \leq \tau \leq 1 } {\mbox{sup}} \| W_d(\tau)\|^2$ is explicitly known.
 In the general case, the quantile of the limit distribution of the first procedure (based on $\widehat{C}_n$) can be computed through Monte-Carlo simulations.
 In the sequel, we will take $q\equiv 1$.
 The Theorem \ref{theo4} below implies that the statistics $\widehat{C}_n $ and $ \widehat{Q}_n $ are too large under the alternative.
 For any $\alpha \in (0,1)$, denote $c_{\alpha}$ the $(1-\alpha)$-quantile of the distribution of  $ \underset{ 0 \leq \tau \leq 1 } {\mbox{sup}} \| W_d(\tau)\|^2$.
 Then at a nominal level $\alpha \in (0,1)$, take $(\widehat{C}_n > c_\alpha )$ as the critical region of the test procedure based on $\widehat{C}_n$.
 This test has correct size  asymptotically.
 On the other hand, it holds that
 \[\underset{ n \rightarrow \infty  } {\limsup }~ P \big( \widehat{Q}_n > c_{\alpha/2} \big) \leq \alpha. \]
 So we can use $c_{\alpha/2}$ as the critical value of the test based on $\widehat{Q}_n$ i.e. $(\widehat{Q}_n > c_\alpha/2 )$ as the critical region.
 This leads to a asymptotically conservative procedure.
 To get correct asymptotic size in the procedure based on $\widehat{Q}_n$, we have to study the asymptotic distribution of $(\widehat{Q}_n^{(1)}, \widehat{Q}_n^{(2)})$.
 This is a very difficult problem due to the dependence structure of the model and the general structure of the parameter.
 In the problem of discriminating between long-range dependence and changes in mean, Berkes {\it et al.} \cite{Berkes2006} have studied the limit
 distribution of such statistic (i.e. the maximum of the maximum between the statistic based on the estimator
 computed with the observations until the time $k$ $(X_1, \ldots,X_k)$ and the one computed with the observations after $k$ $(X_{k+1}, \ldots,X_n)$).
 We have kept this problem as a subject of our future research.

  \subsection{Asymptotic under the alternative}\label{asymp_H1}
Under H$_1$, we assume

\

\noindent {\bf Assumption B }: {\em there exists $\tau_1^*,\ldots,\tau_{K-1}^*$ with $ 0<\tau_1^*<\cdots<\tau_{K-1}^*<1 $ such that for
 $j=1, \ldots,K$, $t_j^{*}=[n\tau_j^*]$ (where $[x]$ is the integer part of $x$).}\\

 \noindent The asymptotic behaviors of these test statistics are given by the following theorem.

\begin{Theo} \label{theo4}
Assume B, $D(\Theta)$, $\mathrm{Id}(\Theta)$,  Var($\Theta$) and $A_i(\Theta)$ $i=0,1,2$ hold with
\[ \sum_{j \geq 1} \sqrt{j} \times \alpha_j^{(i)}(\Theta) < \infty.\]
\begin{enumerate}
  \item Under H$_1$ with $K = 2$, if $\theta_1^* \neq \theta_2^*$ and $\theta_1^*, \theta_2^* \in \ring{\Theta}$
  then \[\widehat{C}_{n} \xrightarrow[n \to +\infty]{P} +\infty .\]
  \item Under H$_1$, if $\theta_1^* \neq \theta_K^*$ and  $\theta_1^*, \theta_2^*, \dots, \theta_K^* \in \ring{\Theta}$ then
  \[\widehat{Q}_n \xrightarrow[n\to +\infty]{P}+\infty . \]
 \end{enumerate}
 \end{Theo}
 It follows that the procedure based on $\widehat{C}_{n}$ is consistent under the alternative of one change while the statistic $\widehat{Q}_n $
  diverges to infinity even if there is multiple under the alternative. So, combined  with an iterated cumulative sums of squares type algorithm
  (see \cite{Tiao1994}) the latter procedure can be used to estimate the number and the break points in the multiple change-points problem.

\begin{figure}
   \begin{center}
   \includegraphics[ width=17.1cm,  height=19.51cm]{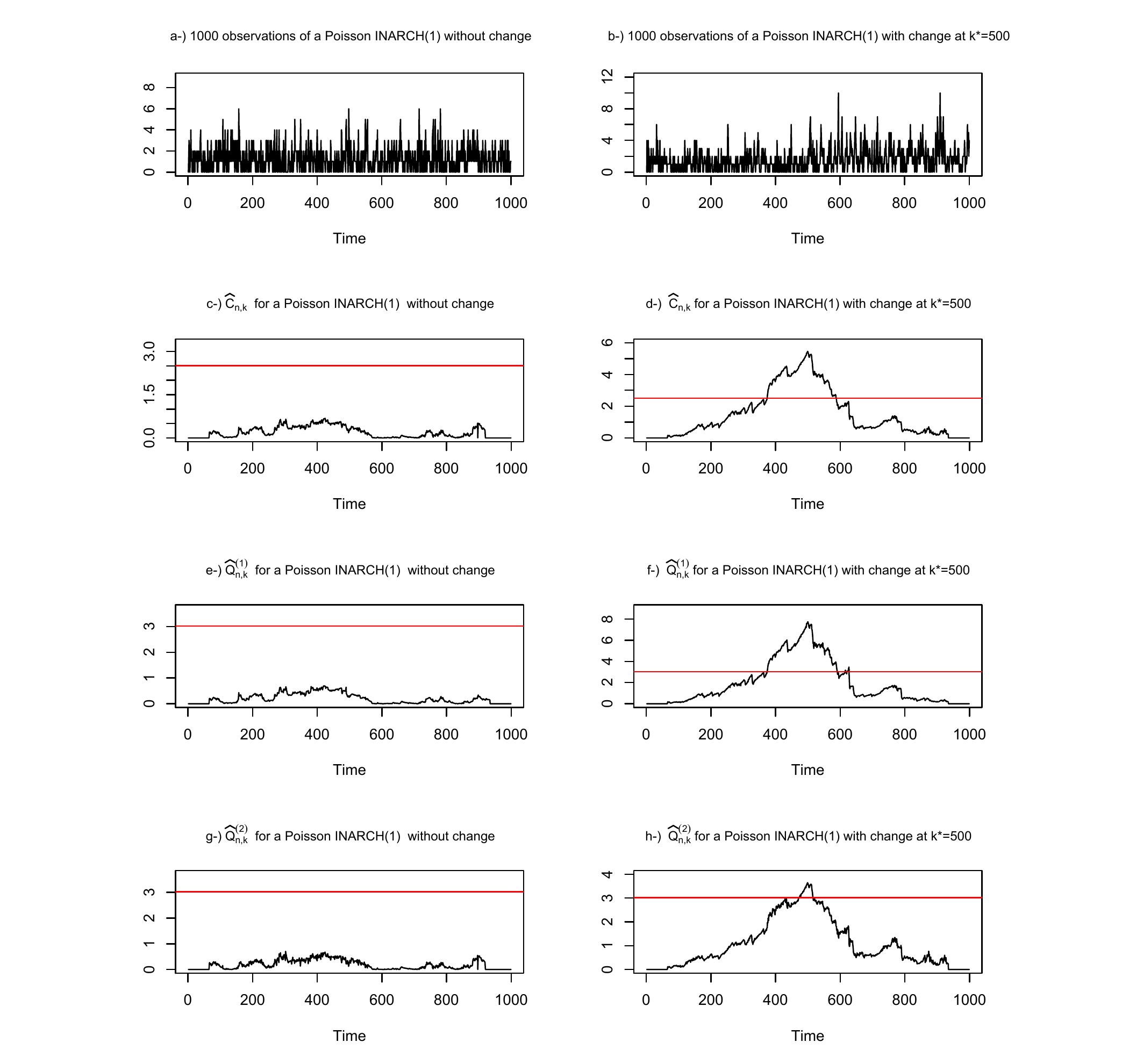}
   \end{center}
   \caption{   Typical realization of 1000 observations of two Poisson INGARCH(1) processes and the corresponding statistics $\widehat{C}^{(1)}_{n,k}$,
    $\widehat{Q}^{(1)}_{n,k}$ and $\widehat{Q}^{(2)}_{n,k}$.
 a-) is a Poisson INGARCH(1) process without change, where the parameter  $\theta_0=(1,0.2)$ is constant.
 b-) is a Poisson INGARCH(1) process where the parameter $\theta_0=(1,0.2)$ changes to $(1,0.45)$ at $k^*=500$.
 c-), d-), e-), f-), g-) and h-)  are their corresponding statistics $\widehat{C}_{n,k}$, $\widehat{Q}^{(1)}_{n,k}$ and $\widehat{Q}^{(2)}_{n,k}$.
 The horizontal line represents the limit of the critical region of the test.}
   \label{fig1}
   \end{figure}

 The Figure \ref{fig1} is an illustration of these tests for an linear Poisson autoregression model of order $1$
 \begin{equation}\label{ill}
   Y_t/\mathcal{F}_{t-1}  \sim  \mbox{ Poisson}(\lambda_t) ~ \text{with} ~ \lambda_t = \alpha_0 +  \beta_1 Y_{t-1}.
   \end{equation}

 \noindent One can see that, under H$_0$ the statistics $\widehat{C}_{n,k}$, $\widehat{Q}^{(1)}_{n,k}$ and $\widehat{Q}^{(2)}_{n,k}$ are below the
 horizontal line (see c-), e-), g-) ) which represents the limit of the critical region. These statistics are greater than the critical value in the neighborhood
 of the breakpoint under H$_1$ (see d-), f-), h-)).
 In several situations, only one of the statistics $\widehat{Q}^{(1)}_{n}$ and $\widehat{Q}^{(2)}_{n}$ is greater than the critical value under the alternative;
 so the use of $\widehat{Q}_n := \max(\widehat{Q}_n^{(1)} , \widehat{Q}_n^{(2)})$ is needed to get more powerful procedure.

 \section{Some simulations results}\label{simu}

  We provide some simulations results to show the empirical performance of these tests procedures.
  We consider a power Poisson INGARCH(1,1) :
  \begin{equation}\label{eq_simu}
   Y_t/\mathcal{F}_{t-1}  \sim  \mbox{ Poisson}(\lambda_t) ~ \text{with} ~ \lambda_t^{\delta} =  \alpha_0 + \alpha_1 \lambda_{t-1}^{\delta}
    +  \beta_1 Y_{t-1}^{\delta}.
   \end{equation}
 For a sample size $n=500, 1000$, the statistics $\widehat{C}_n$ and $\widehat{Q}_n$ are computed with $u_n=v_n=[(\log n)^{5/2}]$.
 The empirical levels and powers reported in the followings table are obtained after 200 replications at the nominal level $\alpha = 0.05$.

  \begin{enumerate}
  \item \textbf{Poisson INGARCH(1,1) with one change-point alternative}. \\
        We assume in (\ref{eq_simu}) that $\delta = 1$ and denote by $\theta = (\alpha_0, \alpha_1, \beta_1)$ the parameter of the model.
        Table \ref{tab1} indicates the empirical levels computed when the parameter is $\theta_0$ and the empirical powers computed when $\theta_0$  changes to
        $\theta_1$ at $n/2$.

   \begin{table}[htbp]
    \centering
    \begin{tabular}{|c c|c|c|c| }
      \hline
           &     & Procedure  &  $n=500$ &  $n=1000$  \\
      \hline
      \hline
   Empirical levels : & $\theta_0=(1,0.1,0.2)$     &  $\widehat{C}_n$ statistic  &  0.020   & 0.040    \\
                            &    & $\widehat{Q}_n$ statistic  &  0.035     &   0.045      \\
                                   &      &    &        &     \\
                      & $\theta_0=(0.3,0.5,0.1)$     &  $\widehat{C}_n$ statistic  &  0.065   &  0.060      \\
                       &      & $\widehat{Q}_n$ statistic  &   0.080    &  0.055       \\
        &    &    &        &     \\
      \hline
     Empirical powers :  & $\theta_0=(1,0.1,0.2)$ ; $\theta_1=(0.7,0.1,0.2)$       & $\widehat{C}_n$ statistic    &  0.415   &  0.840   \\
                          &                 & $\widehat{Q}_n$ statistic    & 0.595    &  0.910    \\
            &      &    &        &       \\
            & $\theta_0=(0.3,0.5,0.1)$ ; $\theta_1=(0.3,0.3,0.4)$       & $\widehat{C}_n$ statistic    &  0.695    &  0.960   \\
                             &                 & $\widehat{Q}_n$  statistic   & 0.865   &   0.995    \\
      \hline
    \end{tabular}
 \caption{{\footnotesize  Empirical levels and powers at the nominal level $0.05$ of test for parameter change in Poisson INGARCH(1,1) model with
                         one change-point alternative.}}
        \label{tab1}
\end{table}

    \item \textbf{Poisson INARCH(1) with two change-points alternative}. \\
          We assume in (\ref{eq_simu}) that $\delta = 1$, $\alpha_1=0$ and denote by $\theta = (\alpha_0, \beta_1)$ the parameter of the model.
          Table \ref{tab2} indicates the empirical levels computed when the parameter is $\theta_0$ and the empirical powers computed when $\theta_0$  changes to
        $\theta_1$ at $0.3n$ which changes to $\theta_2$ at $0.7n$.

   \begin{table}[htbp]
    \centering
    \begin{tabular}{|c c|c|c|c| }
      \hline
           &     & Procedure  &  $n=500$ &  $n=1000$  \\
      \hline
      \hline
   Empirical levels : & $\theta_0 = (1,0.2)$    &  $\widehat{C}_n$ statistic  &  0.040   &  0.055    \\
                            &    & $\widehat{Q}_n$ statistic  &  0.065    & 0.035       \\
                                   &      &    &        &     \\
                      &  $\theta_0 = (0.2,0.5) $    &  $\widehat{C}_n$ statistic  &  0.080    &   0.050     \\
                       &      & $\widehat{Q}_n$ statistic  &  0.090    &    0.050    \\
        &    &    &        &     \\
      \hline
     Empirical powers :  & $\theta_0 = (1,0.2)$ ; $\theta_1 = (1,0.45)$ ; $\theta_2 = (1,0.15)$        & $\widehat{C}_n$ statistic    &  0.715    &   0.955   \\
                          &                 & $\widehat{Q}_n$ statistic    &  0.780    &  0.985   \\
            &      &    &        &       \\
            &  $\theta_0 = (0.2,0.5) $ ; $\theta_1=(0.35,0.5)$ ;  $\theta_2=(0.35,0.4)$     & $\widehat{C}_n$ statistic    &  0.405     &  0.800   \\
                             &                 & $\widehat{Q}_n$  statistic   & 0.540    &   0.865     \\
      \hline
    \end{tabular}
 \caption{{\footnotesize  Empirical levels and powers at the nominal level $0.05$ of test for parameter change in Poisson INARCH(1) model with
                         two change-points alternative.}}
        \label{tab2}
\end{table}

    \item \textbf{ Power Poisson INGARCH(1) with one change-point alternative}. \\
  We assume in (\ref{eq_simu}) that $\delta = 2$, $\alpha_1=0$  and denote by $\theta = (\alpha_0, \beta_1)$ the parameter of the model.
        Table \ref{tab3} indicates the empirical levels computed when the parameter is $\theta_0$ and the empirical powers computed when $\theta_0$  changes to
        $\theta_1$ at $n/2$.

   \begin{table}[htbp]
    \centering
    \begin{tabular}{|c c|c|c|c| }
      \hline
           &     & Procedure  &  $n=500$ &  $n=1000$  \\
      \hline
      \hline
   Empirical levels : & $\theta_0=(0.8,0.15)$     &  $\widehat{C}_n$ statistic  &  0.020    & 0.045   \\
                            &    & $\widehat{Q}_n$ statistic  & 0.065     &   0.055     \\
        &    &    &        &     \\
      \hline
     Empirical powers :  &  $\theta_0=(0.8,0.15)$ ; $\theta_1 = (1.3,0.15)$       & $\widehat{C}_n$ statistic    &  0.415    &  0.870    \\
                          &                 & $\widehat{Q}_n$ statistic    &  0.510    &  0.905     \\
            &      &    &        &       \\
            &  $\theta_0=(0.8,0.15)$ ; $\theta_1 = (0.8,0.4)$     & $\widehat{C}_n$ statistic    & 0.500    &  0.930   \\
                             &                 & $\widehat{Q}_n$  statistic   & 0.590   & 0.945    \\
      \hline
    \end{tabular}
 \caption{{\footnotesize  Empirical levels and powers at the nominal level $0.05$ of test for parameter change in power Poisson INGARCH(1) model with
                         one break alternative.}}
        \label{tab3}
\end{table}
  \end{enumerate}

 It appears in Table \ref{tab1}, \ref{tab2}, \ref{tab3} that these two procedures produces a  size distortion when $n=500$; but
  the empirical levels are close to the nominal one when $n=1000$. One can also see that the empirical powers of these procedures increase with $n$.
  Although the procedure based on $\widehat{Q}_n$ is little more powerful, the test based on $\widehat{C}_n$ provides satisfactory empirical powers even
  in the case of two change-points alternative. This will be the starting point to investigate in our future works, the consistency of this procedure under multiple
  change-points alternative.

 \section{Real data Application}\label{rdata}

\begin{figure}
   \begin{center}
   \includegraphics[width=7.7cm, height=7.51cm]{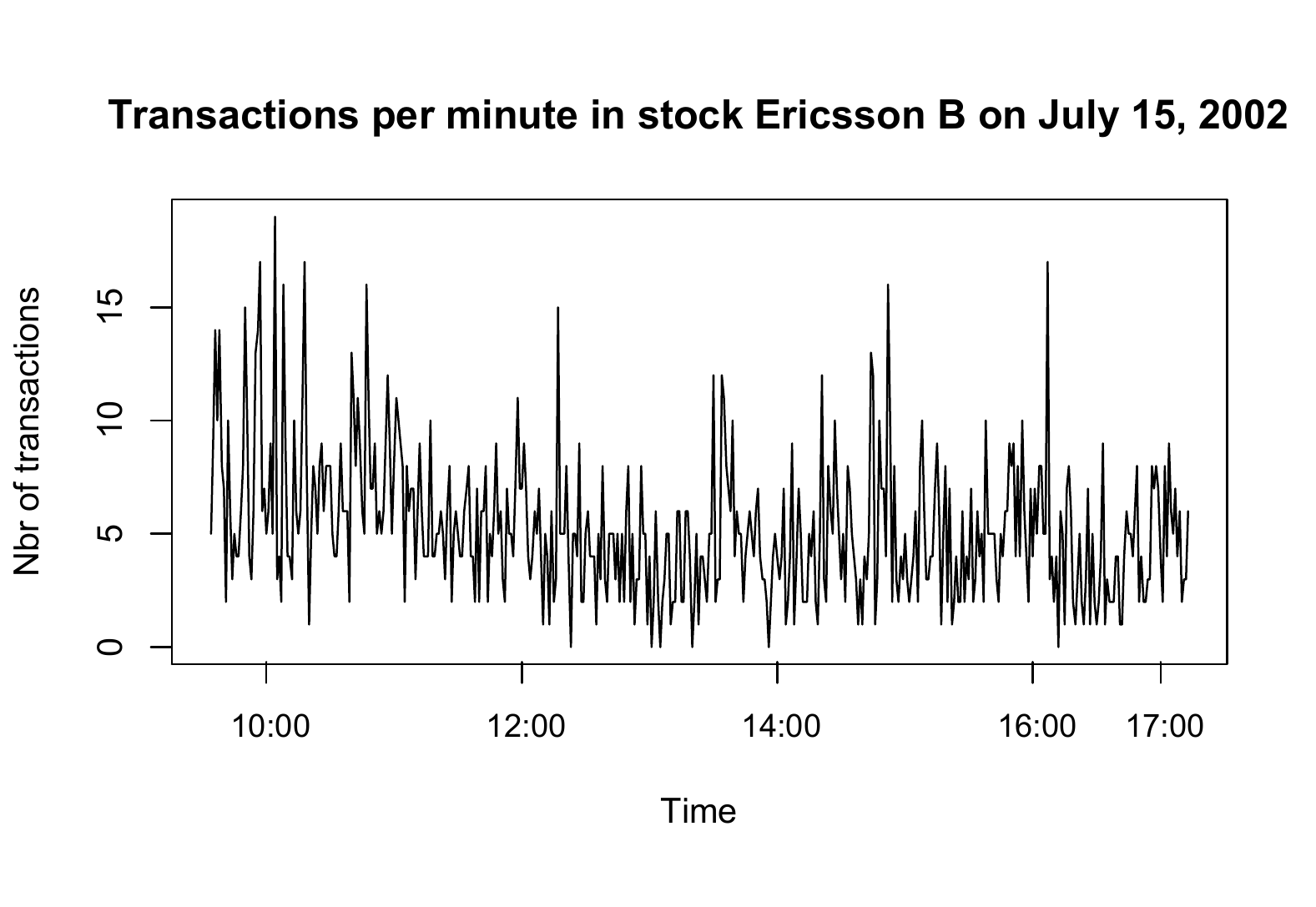}
   \includegraphics[width=7.1cm, height=7.51cm]{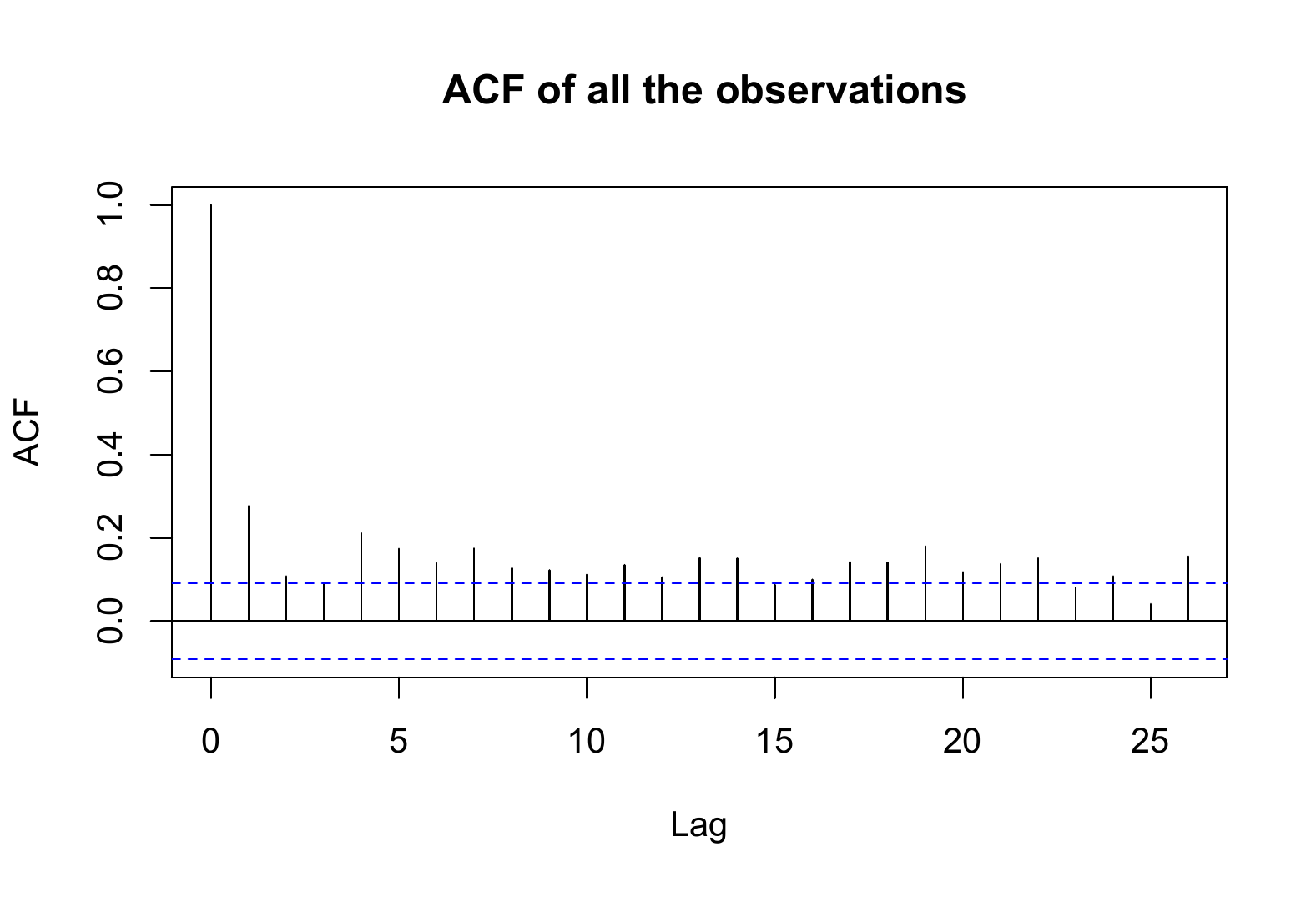}
   \end{center}
   \caption{ Number of transactions per minute for the stock Ericsson B during July 15, 2002 and their autocorrelation function.}
   \label{fig2}
   \end{figure}

   We consider the number of transactions per minute for the stock Ericsson B during July 15, 2002. There are $460$ observations which represent trading from
   09:35 to 17:14.
   Figure \ref{fig2} plot the data and its autocorrelation function.

\begin{figure}
   \begin{center}
    \includegraphics[width=7.7cm, height=5.91cm]{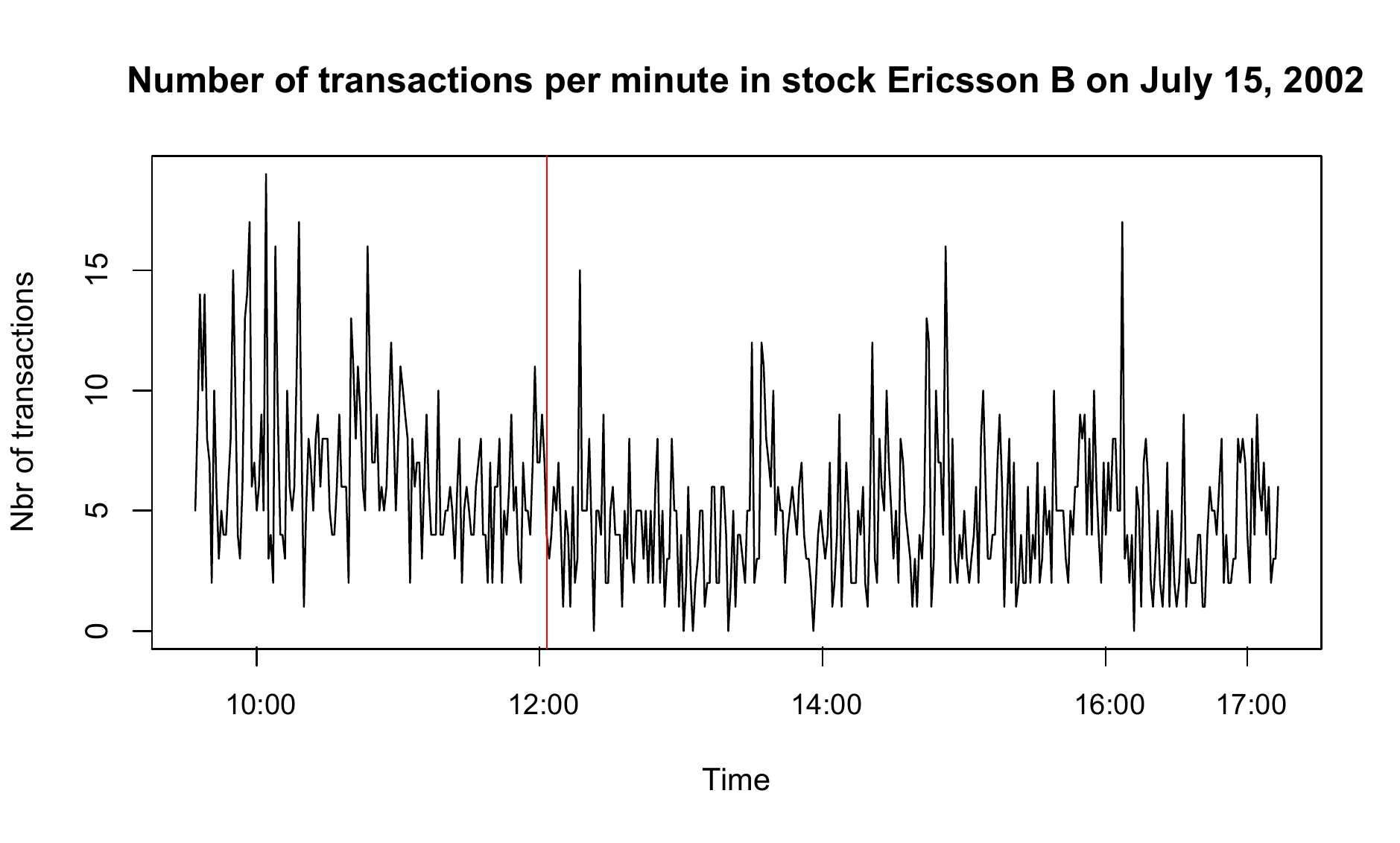}\\
   \includegraphics[width=13.1cm, height=5.31cm]{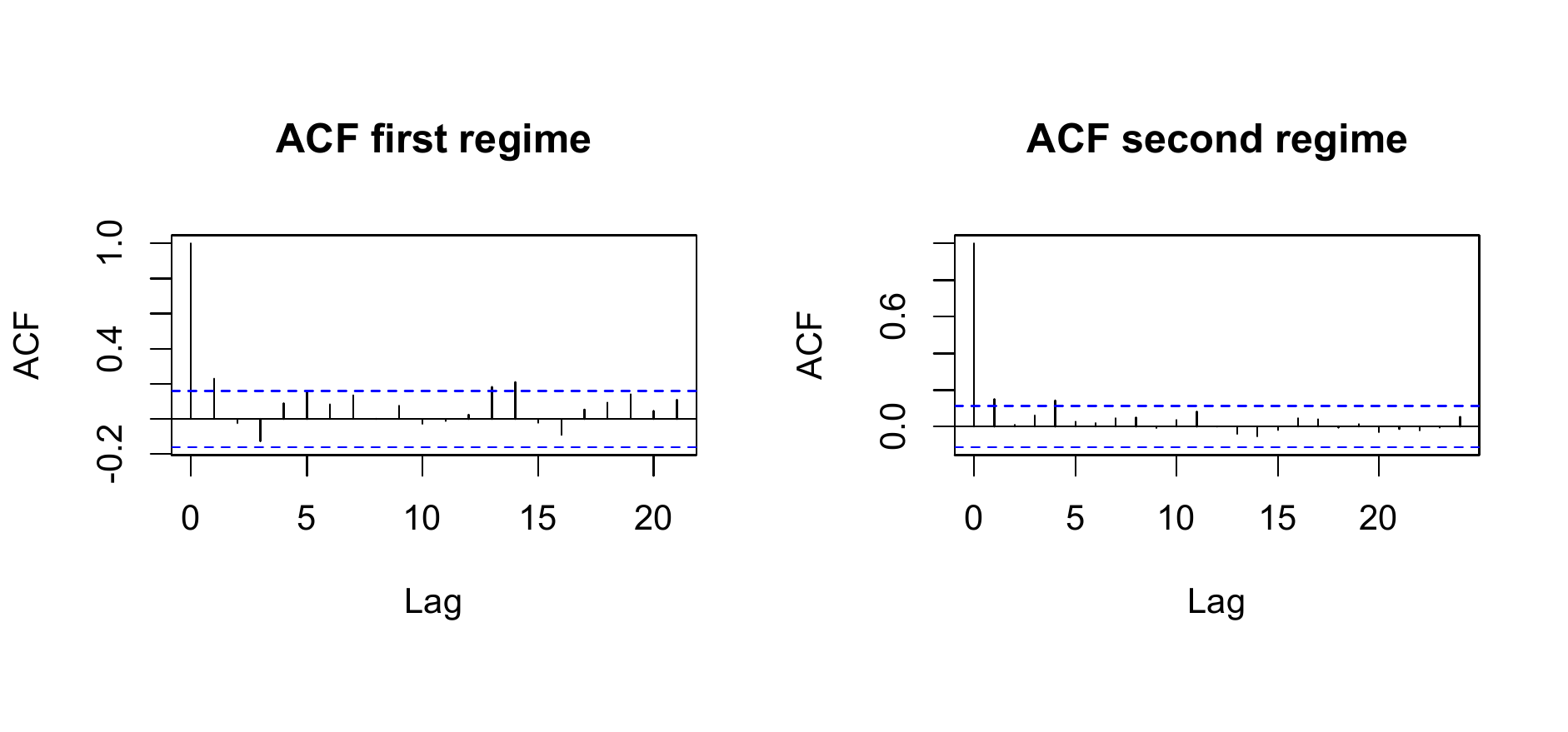}
   \end{center}
   \caption{ Breakpoint in transactions per minute for the stock Ericsson B and the autocorrelation function of each regime.}
   \label{fig3}
   \end{figure}

 Several works (see for instance Fokianos {\it et al.} \cite{Fokianos2009}, Davis and Liu \cite{Davis2012})  on the series of July 2, 2002
 led to use an INGARCH(1,1) model for this series. This model provides $\widehat{\alpha}_1 + \widehat{\beta}_1 $ close to unity. It can be seen in the slow decay
 of the autocorrelation function (see Figure \ref{fig2}). The authors saw a similarity with the high-persistence in IGARCH model.
  The series in the period 2-22 July 2002 have been studied by Brännäs and Quoreshi \cite{Brannas2010}. They have pointed out the presence of long memory
  in these data and applied INARMA model to  both level and first difference forms.  \\
 To test the adequacy of the linearity of the transaction during July 15, 2002,  we have applied the goodness-of-fit test for Poisson count processes
 proposed by Fokianos and Neumann \cite{Fokianos2013}.
 Let $\widehat{\theta}_n = (\widehat{\alpha}_{0,n}, \widehat{\alpha}_{1,n}, \widehat{\beta}_{1,n} )$ be the maximum likelihood estimator computed on the observations.
 Denote $\widehat{I}_t = (\widehat{\lambda}_t, Y_t)$ where
 $\widehat{\lambda}_t = \widehat{\alpha}_{0,n} + \widehat{\alpha}_{1,n}\widehat{\lambda}_{t-1}+ \widehat{\beta}_{1,n}Y_{t-1}$.
 The estimated Pearson residuals is defined by $\widehat{\xi}_t=(Y_t-\widehat{\lambda}_t)/\sqrt{\widehat{\lambda}_t}$. The goodness-of-fit test is based on the
 statistic
 \[  \widehat{T}_n = \sup_{x \in \Pi} |\widehat{G}_n(x)| ~ \text{ with } ~  \widehat{G}_n(x)= \frac{1}{\sqrt{n}} \sum_{t=1}^{n}\widehat{\xi}_t w(x-\widehat{I}_{t-1})\]
 where $\Pi=[0, \infty )^2$ and $w(x) = w(x_1,x_2) = K(x_1)K(x_2)$ where $K(\cdot)$ is a univariate kernel. See \cite{Fokianos2013} for more detail on this test
 procedure.\\
  We have applied this test with uniform and Epanechnikov kernel and the p-values $0.032$ and $0.05$ have been obtained respectively. So, the linear assumption of
  the model is rejected. Recall that, Fokianos and Neumann \cite{Fokianos2013} have already pointed out some doubt about the linearity assumption when they
   analyzed the series of 2 July 2002. \\
 The previous test for change detection have been applied to the series of July 15, 2002. A change has been detected around the midday at $t^*=12:05$.
 Figure \ref{fig3} shows the breakpoint and the autocorrelation function of each regime.
 %
 The previous goodness-of-fit test shows the adequacy of the linear Poisson autoregression of the first regime and raises some doubt about the linearity on the second regime.
 This shows that, the model structure of the transactions in the morning may be different to the structure of the transactions in the afternoon.
 Moreover,  Figure \ref{fig3} shows that the autocorrelation function of each regime decreases fast; this rules out the idea of the long memory in the series.

 \section{Proofs of the main results}\label{proof}
    \noindent {\bf Proof of the Proposition \ref{prop1}} \\
 We will use the same techniques as in \cite{Doukhan2008}.
 Let $p,q$ two fixed non-negative integers. Definite the sequence $(\lambda_t^{p,q})_{t \in \Z}$ by
 \begin{equation}\label{def_lambdapq}
  \lambda_t^{p,q} =
  \begin{cases}
  0 & \text{if } t\leq -q \\
  F(\lambda_{t-1}^{p,q}, \ldots,\lambda_{t-p}^{p,q},0, \ldots;Y_{t-1}, \ldots)  & \text{ otherwise}.
  \end{cases}
  \end{equation}
 The existence of moment of any order of the process $(Y_t,\lambda_t)_{t \in \Z}$ (see \cite{Doukhan2012}) and assumption (\textbf{A}$_F$) imply the
 existence of moment of any order of $(\lambda_t^{p,q})_{t \in \Z}$.
 Let us show that $(\lambda_0^{p,q})_{q\geq0}$ is a Cauchy sequence in $L^1$. By using (\textbf{A}$_F$), we have
 $$ E|\lambda_0^{p,q+1} - \lambda_0^{p,q} | \leq \sum\limits_{j=1}^{p}\alpha_j E|\lambda_{-j}^{p,q+1} - \lambda_{-j}^{p,q} | .$$
 By definition and the strictly stationarity of $(Y_t)_{t \in \Z}$, we can easily see that for $j=1 \ldots,p$, the couples
 $(\lambda_{-j}^{p,q+1},\lambda_{-j}^{p,q+1})$ and $(\lambda_{0}^{p,q+1-j}, \lambda_{0}^{p,q-j})$ have the same distribution. Hence, if comes that
 $$ E|\lambda_0^{p,q+1} - \lambda_0^{p,q} | \leq \sum\limits_{j=1}^{p}\alpha_j E|\lambda_{0}^{p,q+1-j} - \lambda_{0}^{p,q-j} |.$$
  For any fixed $p$, denote $v_q=E|\lambda_0^{p,q+1} - \lambda_0^{p,q}|$ for all $q>p$. It holds that
  $$ v_q \leq \sum\limits_{j=1}^{p}\alpha_j v_{q-j}.$$
  By applying the Lemma 5.4 of \cite{Doukhan2008}, we obtain
  $$v_q\leq \alpha^{q/p}v_0 \text{ where } \alpha = \sum\limits_{j=1}^{\infty}\alpha_j.$$
  Hence, $v_q \rightarrow 0$ as $q \rightarrow \infty$. Thus, for any $p>0$, the sequence $(\lambda_0^{p,q})$ is a Cauchy sequence in $L^1$.
  Therefore, it converges to some limit denoted $\lambda_0^p$. Moreover, since the sequence $(\lambda_0^{p,q})_{q\geq1}$ is measurable w.r.t to
  $\sigma(Y_t, t<0)$, it is the case of the limit $\lambda_0^p$. So, there exists a measurable function $f^{(p)}$ such that
  $\lambda_0^p=f^{(p)}(Y_{-1}, \ldots).$
  By going along similar lines, it holds that for any $t \in \Z$, the sequence $(\lambda_t^{p,q})_{q\geq1}$ converges  in $L^1$ to some
  $\lambda_t^p=f^{(p)}(Y_{t-1}, \ldots)$ and since $(Y_t)_{t \in \Z}$, is stationary and ergodic, the process $(\lambda_t^p)_{t \in \Z}$
  is too stationary and ergodic.\\
  Let $p$ and $t$ fixed. For $q$ large enough, $\lambda_t^{p,q}= F(\lambda_{t-1}^{p,q}, \ldots,\lambda_{t-p}^{p,q},0, \ldots;Y_{t-1}, \ldots)$
  (see (\ref{def_lambdapq})). By using the continuity (which comes from Lipschitz-type conditions)
  of $(Y_1, \ldots,Y_p) \mapsto F(Y_1, \ldots,Y_p,0 \ldots;y)$ for any fixed $y=(y_1)_{i\geq1}$ and by carrying $q$ to infinity, it holds that
  \begin{equation}\label{def_lambdap}
   \lambda_t^{p}= F(\lambda_{t-1}^{p}, \ldots,\lambda_{t-p}^{p},0, \ldots;Y_{t-1}, \ldots).
  \end{equation}
  Denote $\mu_p=E \lambda_t^{p} $, $\mu = \underset{p \geq 1} \sup~ \mu_p$, $\Delta_{p,t}=E|\lambda_t^{p+1} - \lambda_t^{p}|$
  and $\Delta_p = \underset{t \in \Z} \sup \Delta_{p,t}$. By going the same lines as in \cite{Doukhan2008}, we obtain $\Delta_p\leq C \alpha_{p+1}$.
  Therefore, $\Delta_p \rightarrow 0$ as $p\rightarrow \infty$.
  This show that for any fixed $t \in \Z$, $(\lambda_t^{p})_{p\geq1}$ is a Cauchy sequence in $L^1$. Thus it converges to some random $\widetilde{\lambda}_t \in L^1$.
  Moreover, $\widetilde{\lambda}_t$ is measurable w.r.t $\sigma(Y_j, j<t)$ (because it is the case of $(\lambda_t^{p})_{p\geq1}$).
  Thus, there exists a measurable function $f$ such that $\widetilde{\lambda}_t = f(Y_{t-1}, \ldots)$ for any $t \in \Z$.
  This implies that $(\widetilde{\lambda}_t)_{t \in \Z}$ is strictly stationary and ergodic.
  Finally, by using equation (\ref{def_lambdap}) and continuity of $F$, it comes that
  \begin{equation}
   \widetilde{\lambda}_t= F(\widetilde{\lambda}_{t-1}, \ldots;Y_{t-1}, \ldots) \text{ for any } t \in \Z.
  \end{equation}
   Hence, the process $(Y_t,\widetilde{\lambda}_t)_{t \in \Z}$ is strictly stationary ergodic and satisfying (\ref{model1}).
   By the uniqueness of the solution, it holds that $\widetilde{\lambda}_t = \lambda_t $ a.s.
   Thus $\lambda_t=f(Y_{t-1}, \ldots)$ for any $t \in \Z$. \Box

\paragraph{Proof of the Theorem \ref{theo1}}~

Without loss generality, for simplifying notation, we will make the proof with $T_{j_n,k_n} = T_{1,n}$. It will follow two steps.
We will first show that $\norme{\frac{1}{n}\sum_{t\in T_{1,n}}\widehat{\ell}_t(\theta) - L(\theta)}_{\Theta} \xrightarrow[n\to +\infty]{a.s.} 0$
where $L(\Theta):= E(\ell_0(\theta))$; secondly, we will show that the function $\theta \mapsto L(\Theta)$ has a unique maximum in $\theta_0$.
\begin{enumerate}
  \item[(i)] Let $\theta \in \Theta$, recall that $\ell_t(\theta) = Y_t\log\lambda_t(\theta) - \lambda_t(\theta) = Y_t\log f_{\theta}^t - f_{\theta}^t$.
 We have
 \begin{align*}
 \abs{\ell_t(\theta)} &\leq \abs{Y_t}\abs{\log f_{\theta}^t} + \abs{f_{\theta}^t} \\
 &\leq \abs{Y_t}\Big| \log\Big(\frac{f_{\theta}^t}{\underline{c}} \times \underline{c}\Big) \Big | + \abs{f_{\theta}^t} \\
 &\leq \abs{Y_t}\Big(\Big| \frac{f_{\theta}^t}{c}-1\Big| + \abs{\log \underline{c}}\Big) + \abs{f_{\theta}^t} (\esp{for} x > 1, \abs{\log x} \leq \abs{x-1}) \\
 &\leq \abs{Y_t}\left(\frac{1}{\underline{c}}\abs{f_{\theta}^t} + 1 +
 \abs{\log \underline{c}}\right) + \abs{f_{\theta}^t}.
 \end{align*}
 Hence,
 \[\sup_{\theta \in \Theta} \abs{\ell_t(\theta)} \leq \abs{Y_t}\Big(\frac{1}{\underline{c}}\norme{f_{\theta}^t}_{\Theta} + 1 + \abs{\log \underline{c}}\Big)
 + \norme{f_{\theta}^t}_{\Theta}.\]

 We will show that, for any $r > 0$, $\nf(\norme{f_{\theta}^t}^r_{\Theta}) < \infty$. Since A$_0(\Theta)$ holds, we have
 \[ \norme{f^t_{\theta}}_{\Theta} \leq \norme{f^t_{\theta} - f_{\theta}(0)}_{\Theta} + \norme{f_{\theta}(0)}_{\Theta}
 \leq \sum_{j \geq 1} \alpha_j^{(0)}(\Theta)\abs{Y_{t-j}} +
 \norme{f_{\theta}(0)}_{\Theta}.\]
 Thus, by using the stationarity of the process $(Y_t)_{t \in \Z}$, it follows that
 \[\Big( \E(\norme{f^t_{\theta}}_{\Theta}^r) \Big)^{1/r} \leq \norme{Y_0}^r \sum_{j \geq 1} \alpha_j^{(0)}(\Theta) + \norme{f_{\theta}(0)}_{\Theta} < \infty.\]
Therefore, we have
 \[ \nf\Big(\sup_{\theta\in \Theta}\abs{\ell_t(\theta)}\Big) \leq \frac{1}{\underline{c}}(\nf\abs{Y_t}^2)^{1/2}\cdot\left(\nf
 \norme{\ftt}^2_{\Theta}\right)^{1/2} + (1 + \abs{\log \underline{c}})\nf\abs{Y_t} + \nf\norme{\ftt}_{\Theta} < \infty. \]
By the uniform strong law of large number applied on $(\ell_t(\theta))_{t\geq 1}$, it holds that
\begin{equation}\label{ulln_lt}
\Big \| \frac{1}{n}\sum_{t\in T_{1,n}}\ell_t(\theta)-\nf \ell_0(\theta) \Big \|_{\Theta} \xrightarrow[n \to +\infty]{a.s.} 0.
\end{equation}
Now let us show that
\[
\frac{1}{n} \Big \|\sum_{t \in T_{1,n}}\ell_t(\theta)-\sum_{t \in T_{1,n}}\widehat{\ell}_t(\theta)\Big \|_{\Theta}\xrightarrow[n\to +\infty]{a.s.} 0.
 \]
We have
\[\frac{1}{n} \Big \|\sum_{t \in T_{1,n}}\ell_t(\theta)-\sum_{t \in T_{1,n}}\widehat{\ell}_t(\theta) \Big \|_{\Theta} \leq
\frac{1}{n}\sum_{t \in T_{1,n}}\norme{ \ell_t(\theta) - \widehat{\ell}_t(\theta)}_{\Theta}.\]
We will apply the Corollary~$1$ of Kounias and Weng ($1969$). So, it suffices to show that
\[
\frac{1}{n}\sum_{t \geq 1}\frac{1}{t}\nf \Big(\norme{ \ell_t(\theta) - \widehat{\ell}_t(\theta)}_{\Theta}\Big) < \infty.
\]
For $t \in T_{1,n}$ and $\theta \in \Theta$, we have
$\ell_t(\theta)-\widehat{\ell}_t(\theta) = Y_t\log \ftt - \ftt - Y_t\log
\widehat{f}_{\theta}^t + \widehat{f}_{\theta}^t = Y_t(\log f_{\theta}^t-\log
\widehat{f}_{\theta}^t) - (\ftt - \widehat{f}_{\theta}^t)$.

By using the relation $\abs{\log \ftt - \fttc} \leq \frac{1}{\nc}\abs{\ftt - \fttc}$,
it comes that
\[\norme{\ell_t(\theta) - \widehat{\ell}_t(\theta)}_{\Theta} \leq \frac{1}{\nc} \abs{Y_t}\norme{\ftt-\fttc}_{\Theta} + \norme{\ftt-\fttc}_{\Theta}
\leq \Big(\frac{1}{\nc}\abs{Y_t}+1 \Big)\norme{\ftt-\fttc}.\]
By Cauchy-Schwartz inequality,
\begin{align*}
\nf\Big(\norme{\ell_t(\theta)-\widehat{\ell}_t(\theta)}_{\Theta}\Big) &\leq \nf \Big[\Big(\frac{1}{\nc}\abs{Y_t}+1\Big)\norme{\ftt-\fttc}_{\Theta}\Big] \\
&\leq \Big(\nf\Big(\frac{1}{\nc}\abs{Y_t}+1\Big)^{2}\Big)^{1/2}\times \Big(\nf\norme{\ftt-\fttc}^{2}_{\Theta}\Big)^{1/2}.
\end{align*}
We have (by Minkowski inequality),
\[\Big(\nf\Big(\frac{1}{\underline{c}}\abs{Y_t}+1\Big)^{2}\Big)^{1/2} \leq \frac{1}{\underline{c}}(\nf\abs{Y_t}^{2})^{1/2}+1 < \infty.\]
Thus, it comes that
\[\nf \Big(\norme{\ell_t(\theta)-\widehat{\ell}_t(\theta)}_{\Theta}\Big) \leq C \Big(\nf\norme{\ftt-\fttc}^{2}_{\Theta}\Big)^{1/2}.\]
But, we have $\norme{\ftt-\fttc}_{\Theta} \leq \sum_{j \geq
t}\ajo(\Theta)\abs{Y_{t-j}}$.
By using Minkowski inequality, it comes that
\[\Big(\nf \Big \| \ftt-\fttc \Big \|^2_{\Theta} \Big)^{1/2} \leq (\E \abs{Y_0}^2)^{1/2} \sum_{j \geq t}\ajo(\Theta).\]
Hence
\[\Big(\nf\norme{\ftt-\fttc}^2_{\Theta}\Big)^{1/2} \leq C\sum_{j \geq 0}\ajo(\Theta).\]
Thus
\[\nf \Big(\norme{\ell_t(\theta)-\widehat{\ell}_t(\theta)}_{\Theta}\Big) \leq C\sum_{j\geq t}\ajo(\Theta).\]
Therefore
\begin{align*}
\sum_{t \geq 1}\frac{1}{t}\nf\norme{\ltt-\lttc}_{\Theta} &\leq C\sum_{t\geq 1}\frac{1}{t}\sum_{j\geq t}\ajo(\Theta) = C\sum_{t\geq 1}\sum_{j\geq t}\frac{1}{t}\ajo(\Theta) \\
&\leq \sum_{j\geq 1}\sum_{t=1}^j \frac{1}{t}\ajo(\Theta) = C\sum_{j\geq 1}\ajo(\Theta)\sum_{t=1}^j\frac{1}{t} \\
&\leq C\sum_{j=1}\ajo(\Theta)\cdot(1+\log j) \\
&\leq C\sum_{j\leq 1}\ajo(\Theta) + C\sum_{j\geq 1}\ajo(\Theta)\log j \\
&\leq 2C\sum_{j \geq 1}\ajo(\Theta)\log j < \infty .
\end{align*}
Hence, it follows that
\begin{equation}\label{approY_lt}
\frac{1}{n} \Big \|\sum_{t\in T_{1,n}}\ltt - \sum_{t\in T_{1,n}}\lttc \Big \|_{\Theta} \xrightarrow[n\to +\infty]{a.s.} 0.
\end{equation}
From~(\ref{ulln_lt})  and (\ref{approY_lt}),  we deduce that
\[
\Big \|\frac{1}{n}\sum_{t\in T_{1,n}}\lttc - \nf \ell_0(\theta)\Big \|_{\Theta}\xrightarrow[n \to \infty]{a.s.}0.
\]
  \item[(ii)] We  will now show that the function $\theta \mapsto L(\theta) = E\ell_0(\theta)$ has a unique maximum at $\theta_0$. We will
  proceed as in \cite{Davis2012}. Let $\theta \in \Theta$, with $\theta \neq \theta_0$. We have
  \begin{align*}
  L(\theta_0) - L(\theta) &= \nf \ell_0(\theta_0) - \nf \ell_0(\theta) \\
  &= \nf (Y_0\log f_{\theta_0}^0 - f_{\theta_0}^0) - \nf (Y_0\log f_{\theta}^0 - f_{\theta}^0) \\
  &= \nf \left[f_{\theta_0}^0(\log f_{\theta_0}^0-\log f_{\theta_0}^0)\right] - \nf (f_{\theta_0}^0 - f_{\theta_0}^0).
  \end{align*}
 By applying the mean  value theorem at the function $x \mapsto \log x$ defined on $[\underline{c} , +\infty[$, there exists $\xi$ between $f_{\theta_0}^0$ and
 $f_{\theta}^0$ such that
 $\log f_{\theta_0}^0 - \log f_{\theta_0}^0 = ( f_{\theta_0}^0 - f_{\theta_0}^0)\frac{1}{\xi}$.
 Hence,  it comes that
\begin{align*}
L(\theta_0)  - L(\theta) &= \nf \Big(\frac{1}{\xi} f_{\theta_0}^0  (f_{\theta_0}^0-f_{\theta}^0 ) \Big) - \nf (f_{\theta_0}^0 - f_{\theta}^0) \\
&= \nf \Big(\big(\frac{f_{\theta_0}^0}{\xi}-1\big)\big(f_{\theta_0}^0-f_{\theta}^0\big)\Big) \\
&= \nf\Big(\frac{1}{\xi}(f_{\theta_0}^0-\xi)(f_{\theta_0}^0-f_{\theta}^0)\Big).
\end{align*}
Since $\theta \neq \theta_0$, it follows from assumption $\mathrm{Id}(\Theta)$ that $\frac{1}{\xi}(f_{\theta_0}^0-\xi)(f_{\theta_0}^0-f_{\theta}^0) \neq 0$ a.s.
Moreover
\begin{itemize}
  \item if $f_{\theta_0}^0 < f_{\theta}^0$, then $f_{\theta_0}^0 < \xi < f_{\theta}^0$ and hence $\frac{1}{\xi}(f_{\theta_0}^0-\xi)(f_{\theta_0}^0-f_{\theta}^0) > 0$;
  \item if $f_{\theta_0}^0 > f_{\theta}^0$, then $f_{\theta}^0 < \xi < f_{\theta_0}^0$ and hence $\frac{1}{\xi}(f_{\theta_0}^0-\xi)(f_{\theta_0}^0-f_{\theta}^0) > 0$.
\end{itemize}
We deduce that $\frac{1}{\xi}(f_{\theta_0}^0-\xi)(f_{\theta_0}^0-f_{\theta}^0) > 0$ a.s.. Hence
$L(\theta_0) - L(\theta) = \nf \Big(\frac{1}{\xi}(f_{\theta_0}^0-\xi)(f_{\theta_0}^0-f_{\theta}^0)\Big) > 0$. Thus, the
function $\theta \mapsto L(\theta)$ has a unique maximum at $\theta_0$.
\end{enumerate}
(i), (ii) and  standard arguments lead to the consistency of $\widehat{\theta}_n(T_{1,n})$.   \Box

~\\
   The following lemma is needed to prove the Theorem \ref{theo2}.

 \begin{lem}\label{lem1}
 Let $(j_n)_{n\geq 1}$ and $(k_n)_{n\geq 1}$ two integer valued sequences such that $(j_n)_{n\geq 1}$ is increasing, $j_n \rightarrow \infty$ and
 $k_n - j_n \rightarrow \infty$ as $n\rightarrow \infty$. Let $n\geq1$, for any segment $T=T_{j_n,k_n} \subset \{1, \ldots,n \}$, it holds under assumptions
 of Theorem \ref{theo2} that

\begin{itemize}
  \item[(i)] $\E \Big( \dfrac{1}{\sqrt{k_n - j_n}} \Big \| \dfrac{\partial}{\partial \theta} L_n(T,\theta)
  - \dfrac{\partial}{\partial \theta} \widehat{L}_n(T,\theta)\Big \|_{\Theta} \Big)  \underset{n \rightarrow \infty}{\longrightarrow 0}$ ;
  \item[(ii)]  $ \Big \| \dfrac{1}{k_n - j_n} \dfrac{\partial^2}{\partial \theta \partial \theta ' } L_n(T,\theta) -
  \E \Big( \dfrac{\partial^2 \ell_0(\theta) }{\partial \theta  \partial \theta ' } \Big) \Big \|_{\Theta} \underset{n \rightarrow \infty }{\overset{a.s.} \longrightarrow 0}$ ;
  \item[(iii)] $ \Big \|\dfrac{1}{k_n - j_n} \sum_{t \in T}\dfrac{1}{\fttc}\big(\dfrac{\partial}{\partial\theta}\fttc\big)\big(\dfrac{\partial}{\partial\theta}\fttc\big)'
    - \E\Big( \dfrac{1}{f^0_{\theta}}\big(\dfrac{\partial}{\partial\theta} f^0_{\theta} \big)\big(\dfrac{\partial}{\partial\theta} f^0_{\theta} \big)' \Big \|_{\Theta} \underset{n \rightarrow \infty }{\overset{a.s.} \longrightarrow 0}$.
  \end{itemize}

\end{lem}

 \begin{dem}

\begin{itemize}
  \item[(i)] Let $i \in \{1, \dots, n\}$. We have
  \[\ddi L_n(T,\theta) = \sum_{t\in T}\ddi (Y_t\log \ftt - \ftt) = \sum_{t\in T}\left(Y_t\frac{1}{\ftt}\ddi\ftt - \ddi\ftt\right)
   = \sum_{t\in T}\ddi \ell_t(\theta)\]
   and
   \[\ddi \widehat{L}_n(T,\theta) = \sum_{t\in T}\left(Y_t\frac{1}{\fttc}\ddi\fttc - \ddi\fttc\right) = \sum_{t\in T}\ddi \widehat{\ell}_t(\theta).\]
   Hence
   \begin{align}
  \nonumber \abs{\ddi \ell_t(\theta) - \ddi \widehat{\ell}_t(\theta)} &= \abs{Y_t\frac{1}{\ftt}\ddi\ftt - \ddi\ftt - Y_t\frac{1}{\fttc}\ddi\fttc + \ddi\fttc}\\
   &\leq \abs{Y_t}\abs{\frac{1}{\ftt}\ddi\ftt - \frac{1}{\fttc}\ddi\fttc} + \abs{\ddi\ftt-\ddi\fttc}. \label{approY_dlt}
   \end{align}
 Using the relation $\abs{a_1b_1 - a_2b_2} \leq \abs{a_1-a_2}\abs{b_2} + \abs{b_1-b_2}\abs{a_1}$, we have
 \begin{align*}
 \abs{\frac{1}{\ftt}\ddi\ftt - \frac{1}{\fttc}\ddi\fttc} &\leq \abs{\frac{1}{\ftt}-\frac{1}{\fttc}}\abs{\ddi\fttc} + \abs{\ddi\ftt-\ddi\fttc}\abs{\frac{1}{\ftt}} \\
 &\leq \frac{1}{\nc^2}\abs{\ftt-\fttc}\abs{\ddi\fttc} + \frac{1}{\nc}\abs{\ddi\ftt - \ddi\fttc}.
 \end{align*}
 Hence, (\ref{approY_dlt}) implies
 \begin{align*}
 \norme{\ddi \ell_t(\theta) - \ddi\widehat{\ell}_t(\theta)} &\leq
 \abs{Y_t}\left(\frac{1}{\nc^2}\norme{\ftt-\fttc}_{\Theta}\norme{\ddi\fttc}_{\Theta}
 + \frac{1}{\nc}\norme{\ddi\ftt - \ddi\fttc}_{\Theta}\right) + \norme{\ddi\ftt - \ddi\fttc}_{\Theta} \\
 &\leq C\abs{Y_t}\norme{\ddi\fttc}_{\Theta} \norme{\ftt - \fttc}_{\Theta} + C(1+\abs{Y_t})\norme{\ddi\ftt - \ddi\fttc}_{\Theta}.
 \end{align*}
 Let $r \gr 0$. Using the Minkowski and Hölder inequalities, it holds that
 \begin{align*}
 \left(E\left[\norme{\ddi \ell_t(\theta) - \ddi \widehat{\ell}_t(\theta)}^r_{\Theta}\right]\right)^{1/r} & \leq
 C\left(E\left[\abs{Y_t}^r\norme{\ddi\fttc}^r_{\Theta}\norme{\ftt-\fttc}^r_{\Theta}\right]\right)^{1/r} \\
 &\qquad + C\left(E\left[(1+\abs{Y_t})^r\norme{\ddi\ftt - \ddi\fttc}^r_{\Theta}\right]\right)^{1/r} \\
 &\leq
 C\left((E\abs{Y_t}^{3r})^{1/3}\left(E\norme{\ddi\fttc}^{3r}_{\Theta}\right)^{1/3}
 \left(E\norme{\ftt-\fttc}^{3r}_{\Theta}\right)^{1/3}\right)^{1/r} \\
 &\qquad + C\left(\left(E(1+\abs{Y_t})^{2r}\right)^{1/r}\left(E\norme{\ddi\ftt - \ddi\fttc}^{2r}\right)^{1/2}\right)^{1/r} \\
 &\leq C\norme{Y_t}_{3r}\norme{\norme{\ddi\fttc}_{\Theta}}_{3r}\cdot\left(E\norme{\ftt-\fttc}_{\Theta}^{3r}\right)^{1/3r} \\
 &\qquad  + C\norme{1+\abs{Y_t}}_{2r}\cdot\left(E\norme{\ddi\ftt -
 \ddi\fttc}_{\Theta}^{2r}\right)^{1/2r}.
 \end{align*}
 But we have $\norme{Y_t}_{3r} = C \pet \infty$ and $\norme{1 + \abs{Y_t}}_{1r} \pet \infty$. Hence
\[ \norme{\ddi\fttc}_{\Theta} \leq \norme{\ddi f_{\theta}(0)}_{\Theta} + \norme{\ddi\fttc-\ddi f_{\theta}(0)}_{\Theta}
 \leq C + \sum_{j \geq 1}\alpha_j^{(1)}(\Theta)\abs{Y_{t-j}}. \]
Thus, \[\norme{\norme{\ddi\fttc}_{\Theta}}_{3r} \leq C +
\norme{Y_0}_{3r}\sum_{j\geq 1}\alpha_j^{(1)}(\Theta) \leq C(1 +
\sum_{j\geq 1}\alpha_j^{(1)}(\Theta)) \pet \infty.\]
 We also have $ \norme{\ftt-\fttc}_{\Theta} \leq \sum_{j \geq t}\alpha_j^{(0)}(\Theta)\abs{Y_{t-j}}$.
 Hence
\[\left(E\norme{\ftt-\fttc}^{3r}_{\Theta}\right)^{1/3r} = \norme{\norme{\ftt-\fttc}_{\Theta}}_{3r} \leq C\sum_{j\geq t}\alpha_j^{(0)}(\Theta).\]
The same argument gives
\[\left(E\norme{\ddi\ftt - \ddi\fttc}^{2r}_{\Theta}\right)^{1/2r} = \norme{\norme{\ddi\ftt - \ddi\fttc}_{\Theta}}_{2r} \leq
C\sum_{j\geq t}\alpha_j^{(1)}(\Theta).\]
Hence,
\[\left(E\left[\norme{\ddi \ell_t(\theta)-\ddi\widehat{\ell}_t(\theta)}^{r}_{\Theta}\right]\right)^{1/r} \leq C\sum_{j\geq t}\left(\alpha_j^{(0)}(\Theta)+\alpha_j^{(1)}(\Theta)\right).\]
Therefore, we have (with $r = 1$)
\begin{align*}
& E\left(\frac{1}{\sqrt{k_n-j_n}}\norme{\ddi L_n(T,\theta) - \ddi \widehat{L}_n(T,\theta)}_{\Theta}\right)
 \leq \frac{1}{\sqrt{k_n-j_n}} \sum_{t\in T}E\norme{\ell_t(\theta)-\widehat{\ell}_t(\theta)}_{\Theta} \\
& \quad \leq C\frac{1}{\sqrt{k_n-j_n}} \sum_{t\in T}\left(\sum_{j\geq t}\left(\alpha_j^{(0)}(\Theta)+\alpha_j^{(1)}(\Theta)\right)\right)\\
& \quad\leq C\frac{1}{\sqrt{k_n-j_n}} \sum_{t\in T} \left[\sum_{j=t}^{k_n}\left(\alpha_j^{(0)}(\Theta)+\alpha_j^{(1)}(\Theta)\right) + \sum_{j\geq k_n}\left(\alpha_j^{(0)}(\Theta)+\alpha_j^{(1)}(\Theta)\right)\right]\\
& \quad\leq C\frac{1}{\sqrt{k_n-j_n}} \left[\sum_{t=j_n}^{k_n}\sum_{j=t}^{k_n}\left(\alpha_j^{(0)}(\Theta)+\alpha_j^{(1)}(\Theta)\right) + \sum_{t=j_n}^{k_n}\sum_{j\geq k_n}\left(\alpha_j^{(0)}(\Theta)+\alpha_j^{(1)}(\Theta)\right)\right] \\
& \quad\leq C\frac{1}{\sqrt{k_n-j_n}} \left[\sum_{j=j_n}^{k_n}\sum_{t=j_n}^{j}\left(\alpha_j^{(0)}(\Theta)+\alpha_j^{(1)}(\Theta)\right) + (k_n-j_n)\sum_{j\geq k_n}\left(\alpha_j^{(0)}(\Theta)+\alpha_j^{(1)}(\Theta)\right)\right]\\
& \quad \leq C\frac{1}{\sqrt{k_n-j_n}} \left[\sum_{j=j_n}^{k_n}(j-j_n)\left(\alpha_j^{(0)}(\Theta)+\alpha_j^{(1)}(\Theta)\right) + (k_n-j_n)\sum_{j\geq k_n}\left(\alpha_j^{(0)}(\Theta)+\alpha_j^{(1)}(\Theta)\right)\right] \\
& \quad \leq C\frac{1}{\sqrt{k_n-j_n}} \sum_{j=j_n}^{k_n}(j-j_n)\left(\alpha_j^{(0)}(\Theta)+\alpha_j^{(1)}(\Theta)\right) + C\sqrt{k_n-j_n}\sum_{j\geq k_n}\left(\alpha_j^{(0)}(\Theta)+\alpha_j^{(1)}(\Theta)\right)\\
& \quad \leq \frac{C}{\sqrt{k_n-j_n}}  \sum_{j=j_n}^{j_n+\log(k_n-j_n)}(j-j_n)\left(\alpha_j^{(0)}(\Theta)+\alpha_j^{(1)}(\Theta)\right)\\
& \qquad + \frac{C}{\sqrt{k_n-j_n}} \sum_{j = j_n+\log(k_n-j_n)}^{k_n}(j-j_n)\left(\alpha_j^{(0)}(\Theta)+\alpha_j^{(1)}(\Theta)\right)
+ C\sum_{j\geq k_n}\sqrt{j}\left(\alpha_j^{(0)}(\Theta)+\alpha_j^{(1)}(\Theta)\right)
\end{align*}
\begin{align*}
& \quad \leq \frac{C\log(k_n-j_n)}{\sqrt{k_n-j_n}}  \sum_{j\geq 1}\left(\alpha_j^{(0)}(\Theta)+\alpha_j^{(1)}(\Theta)\right) + C\sum_{j \geq j_n+\log(k_n-j_n)}\sqrt{j-j_n}\left(\alpha_j^{(0)}(\Theta)+\alpha_j^{(1)}(\Theta)\right) \\
& \qquad + C\sum_{j\geq k_n}\sqrt{j}\left(\alpha_j^{(0)}(\Theta)+\alpha_j^{(1)}(\Theta)\right) \\
& \qquad \qquad \qquad \xrightarrow[n\to +\infty]{}0.
\end{align*}
This holds for any coordinate $i = 1\dots,d$. This (i) has been proved.

  \item[(ii)] For $i,j \in \{1, \ldots,d \}$, we have
  \begin{align}
\nonumber \dcij \ell_t(\theta) &= \ddj(Y_t\frac{1}{\ftt}\ddi\ftt-\ddi\ftt) \\
\nonumber &= Y_t\left[\ddj\left(\frac{1}{\ftt}\right)\times \ddi\ftt + \frac{1}{\ftt}\times \dcij \ftt\right]-\dcij\ftt \\
\nonumber &= Y_t\left[-\frac{1}{\left(\ftt\right)^2}\left(\ddj\ftt\right)\times \ddi\ftt + \frac{1}{\ftt}\times\dcij\ftt\right]-\dcij\ftt \\
&=
-\frac{Y_t}{\left(\ftt\right)^2}\left(\ddi\ftt\right)\left(\ddj\ftt\right)
+ \left(\frac{Y_t}{\ftt}-1\right)\dcij\ftt. \label{d2_lt}
\end{align}
We will show that $E\left[\norme{\dcij \ell_t(\theta)}\right] \pet +\infty$.
We have
\begin{align*}
E\left[\norme{\dcij \ell_t(\theta)}_{\Theta}\right] &\leq \frac{1}{\nc^2}E\left[\abs{Y_t}\norme{\ddi\ftt}_{\Theta}\norme{\ddj\ftt}_{\Theta}\right] + E\left[\left(\frac{\abs{Y_t}}{\nc}+1\right)\norme{\dcij\ftt}_{\Theta}\right]\\
&\leq
C\norme{Y_t}_3\norme{\norme{\ddi\ftt}_{\Theta}\norme{\ddj\ftt}_{\Theta}}_3
+ C(\norme{Y_t}_2+1)\norme{\norme{\dcij\ftt}_{\Theta}}_2.
\end{align*}
But, we have $\norme{Y_t}_3 = \norme{Y_0}_3 \pet \infty$, and
$\norme{Y_t}_2 \pet \infty$.
\begin{align*}
\norme{\norme{\ddi\ftt}_{\theta}}_3 &\leq \norme{\norme{\ddi f_{\theta}(0)}_{\Theta}}_3 + \norme{\norme{\ddi\ftt - \ddi f_{\theta}(0)}_{\Theta}}_3 \\
&\leq \norme{\ddi f_{\theta}(0)}_{\Theta} + \norme{\sum_{j \geq 1}\alpha_j^{(1)}(\Theta)\abs{Y_{t-j}}}_3 \\
&\leq \norme{\ddi\ftt}_{\Theta} + \norme{Y_0}_3\sum_{j \geq
1}\alpha_j^{(1)}(\Theta) \pet +\infty.
\end{align*}
Similarly, we have $\norme{\norme{\ddj\ftt}_{\Theta}}_3 \pet +\infty$. Using the same argument, we obtain
\[\norme{\norme{\dcij\ftt}_{\Theta}}_2 \leq \norme{Y_0}_2\sum_{j \geq 1}\alpha_j^{(2)}(\Theta) \pet +\infty. \]
Hence, $E\left[\norme{\dcij \ell_t(\theta)}_{\Theta}\right] \pet +\infty$. Thus, for the stationary ergodicity of the sequence
$\left(\dcij \ell_t(\theta)\right)_{t\in \Z}$ and the uniform strong law of large numbers, it holds that
\[\norme{\frac{1}{k_n-j_n}\dcij \ell_t(\theta)- E\dcij \ell_0(\theta)}_{\Theta} \xrightarrow[n \to +\infty]{a.s.} 0.\]
This completes the proof of (ii).

  \item[(iii)]  Goes the same lines as in (i) and (ii).
 \end{itemize}

 \end{dem}

 \paragraph{Proof of Theorem~\ref{theo2}}~
 Here again, without loss generality we will make the proof with $T_{j_n,k_n}=T_{1,n}$.
 Recall that $\Theta \subset \R^d$. Let $T \subset \{1, \ldots,n \}$; for any $\overline{\theta} \in \Theta$ and $i=1, \ldots,n$, by applying the Taylor expansion to
 the function $\theta \mapsto \ddt L_n(T,\theta)$, there exists $\theta_{n,i}$ between $\overline{\theta}$ and $\theta_0$ such that
 \[ \ddi L_n(T,\overline{\theta}) = \ddi L_n(T,\theta_0) + \frac{\partial^2}{\partial\theta \partial\theta_i} L_n(T,\theta_{n,i})\cdot (\overline{\theta} - \theta_0) .\]
 Denote
 \[G_n(T,\overline{\theta}) = - \frac{1}{Card(T)} \Big( \frac{\partial^2}{\partial\theta \partial\theta_i} L_n(T,\theta_{n,i}) \Big)_{1\leq i \leq d}. \]
 It comes that
 \begin{equation} \label{Taylor_G_1}
 Card(T)G_n(T,\overline{\theta})\cdot (\overline{\theta} - \theta_0) = \ddt L_n(T,\theta_0) - \ddt L_n(T,\overline{\theta}) .
 \end{equation}
  By applying (\ref{Taylor_G_1}) with $\overline{\theta} = \widehat{\theta}_n(T)$ we obtain
 \begin{equation} \label{Taylor_G_2}
 Card(T) G_n(T,\widehat{\theta}_n(T))\cdot (\widehat{\theta}_n(T) - \theta_0) = \ddt L_n(T,\theta_0) - \ddt L_n(T,\widehat{\theta}_n(T)) .
 \end{equation}
  (\ref{Taylor_G_2})  holds for any $T \subset \{1, \ldots,n \}$, thus
 \begin{equation} \label{Taylor_G_3}
 \sqrt{n}G_n(T_{1,n},\widehat{\theta}_n(T_{1,n}))\cdot (\widehat{\theta}_n(T_{1,n}) - \theta_0)
 = \frac{1}{\sqrt{n}}\Big(\ddt L_n(T_{1,n},\theta_0) - \ddt L_n(T_{1,n},\widehat{\theta}_n(T_{1,n})) \Big) .
 \end{equation}
  We can rewrite (\ref{Taylor_G_3}) as
\begin{align*}
\sqrt{n}G_n(T_{1,n},\tnc(T_{1,n}))(\tnc(T_{1,n})-\theta_0) &= \frac{1}{\sqrt{n}}\ddt L_n(T_{1,n},\theta_0) - \frac{1}{\sqrt{n}}\ddt \widehat{L}_n(T_{1,n},\tnc(T_{1,n})) \\
& \qquad + \frac{1}{\sqrt{n}}\left( \ddt \widehat{L}_n(T_{1,n},\tnc(T_{1,n})) - \ddt L_n(T_{1,n},\tnc(T_{1,n})) \right).
\end{align*}
For $n$ large enough, $\ddt \widehat{L}_n(T_{1,n},\tnc(T_{1,n})) = 0$, because $\tnc(T_{1,n})$ is a local maximum of $\theta \mapsto \widehat{L}(T_{1,n},\theta)$.
Moreover, according to Lemma \ref{lem1} (i), it holds that
\begin{align*}
E\left(\frac{1}{\sqrt{n}}\abs{\ddt L_n(T_{1,n},\tnc(T_{1,n})) - \ddt \widehat{L}_n(T_{1,n},\tnc(T_{1,n}))}\right) &\leq
E\left(\frac{1}{\sqrt{n}} \norme{\ddt L_n(T_{1,n},\theta) - \ddt \widehat{L}_n(T_{1,n},\theta)}_{\Theta} \right) \\
& \qquad \xrightarrow[n \to +\infty]{} 0.
\end{align*}
So, for $n$ large enough, we have
\begin{equation} \label{approx1_G}
\sqrt{n}G_n(T_{1,n},\tnc(T_{1,n}))(\tnc(T_{1,n})-\theta_0) = \frac{1}{\sqrt{n}}\ddt L_n(T_{1,n},\theta_0) + o_P(1).
\end{equation}
To complete the proof of Theorem \ref{theo2}, we have to show that
\begin{enumerate}
  \item[(a)]  $\left(\ddt \ell_t(\theta_0),\mathcal{F}_t\right)_{t \in \Z}$  is a stationary ergodic martingable
  difference sequence and $E\left(\norme{\ddt \ell_t(\theta_0)}^2\right) \pet \infty$;
  \item[(b)]  $ \Sigma = -E\left(\frac{\partial^2}{\partial\theta\partial\theta'}\ell_0(\theta_0)\right) $ and
    $G_n(T_{1,n},\tnc(T_{1,n})) \xrightarrow[n \to +\infty]{a.s.} \Sigma$;
  \item[(c)] $\Sigma = E\Big[\frac{1}{f^0_{\theta_0}}\left(\frac{\partial}{\partial\theta}f^0_{\theta_0}\right)\left(\frac{\partial}{\partial\theta}f^0_{\theta_0}\right)'\Big] $
   is invertible.
\end{enumerate}

\begin{enumerate}
  \item[(a)] Recall that $\ddt \ell_t(\theta_0) = \left(\frac{Y_t}{f_{\theta_0}^t}-1\right)\ddt f_{\theta_0}^t$
  and $\mathcal{F}_t = \sigma(Y_s, s\leq t)$. Since the functions $f_{\theta_0}^t$ and $\ddt f_{\theta_0}^t$ are $\mathcal{F}_{t-1}$-measurable,
   we have
   \[E\Big(\ddt \ell_t(\theta_0)\big| \mathcal{F}_{t-1}\Big) =
  \Big(\frac{1}{f_{\theta_0}^t} E(Y_t | \mathcal{F}_{t-1}) - 1\Big)\ddt f_{\theta_0}^t = 0.\]
  Moreover, since $\abs{Y_t}$ and $\norme{\ddt f_{\theta}^t}$ have moment of any order, we have
  \[E\Big( \big|\ddt \ell_t(\theta_0)\big|^2\Big) \leq E\Big(\big(\frac{\abs{Y_t}}{\nc}+1\big)^2\big \| \ddt \ftt \big \|^2_{\Theta}\Big) \pet \infty.\]
  \item[(b)]

  According to (\ref{d2_lt}), we have
  \begin{equation} \label{d2_lt2}
     \dfrac{\partial^2}{\partial \theta  \partial \theta' } \ell_t(\theta) = -\dfrac{Y_t}{ \Big(\ftt\Big)^2} \Big(\ddt\ftt\Big) \Big(\ddt\ftt\Big)'
  + \Big(\dfrac{Y_t}{\ftt}-1\Big)  \dfrac{\partial^2}{\partial \theta  \partial \theta' }  \ftt .
  \end{equation}
  But by using the same argument as in (a), we obtain
  $ \E \Big( \big(\dfrac{Y_t}{ f^0_{\theta_0}}-1\big)  \dfrac{\partial^2}{\partial \theta  \partial \theta' }  f^0_{\theta_0} \big | \mathcal{F}_{t-1}\Big) = 0$.
  Hence, (\ref{d2_lt2}) implies
  \[  \E\Big( \dfrac{\partial^2}{\partial \theta  \partial \theta' } \ell_t(\theta_0) \Big) = - \E \left( \dfrac{Y_t}{(f^0_{\theta_0} )^2}\Big(\ddt f^0_{\theta_0} \Big)\Big(\ddt f^0_{\theta_0}  \Big)' \right)
  = -\E \left( \dfrac{1}{f^0_{\theta_0}} \Big(\ddt f^0_{\theta_0} \Big)\Big(\ddt f^0_{\theta_0}  \Big)' \right)= -\Sigma  .\]
  Moreover, recall that
\[
  G_n(T_{1,n},\tnc(T_{1,n})) = -\frac{1}{n}\left(\frac{\partial^2}{\partial\theta\partial\theta_i} L(T_{1,n},\theta_{n,i})\right)_{1 \leq i \leq d} \\
  = -\frac{1}{n}\left(\sum_{t=1}^n \frac{\partial^2}{\partial\theta\partial\theta_i} \ell_t(\theta_{n,i})\right)_{1 \leq i \leq d}.
\]
  For any $j = 1,\dots,d$, we have
\begin{align*}
\abs{\frac{1}{n}\sum_{t=1}^n\dcji \ell_t(\theta_{n,i}) - E\left(\dcji \ell_0(\theta_0)\right)} &\leq
\abs{\frac{1}{n}\sum_{t=1}^n\dcji \ell_t(\theta_{n,i}) - E\left(\dcji \ell_0(\theta_{n,i})\right)} \\
&\qquad + \abs{E\left(\dcji \ell_0(\theta_{n,i})\right) - E\left(\dcji \ell_0(\theta_0)\right)} \\
&\leq \abs{E\left(\dcji \ell_0(\theta_{n,i})\right) - E\left(\dcji \ell_0(\theta_0)\right)} \\
&\qquad + \norme{\frac{1}{n}\sum_{t=1}^n\dcji \ell_t(\theta) - E\left(\dcji \ell_0(\theta)\right)}_{\Theta} \\
& \quad \xrightarrow[n \to +\infty]{a.s.} 0.
\end{align*}
This holds for any $1 \leq i,j\leq d$. Thus,
\begin{align*}
G_n(T_{1,n},\tnc(T_{1,n})) &= -\frac{1}{n}\left(\frac{\partial^2}{\partial\theta\partial\theta_i}\ell_t(\theta_{n,i})\right)_{1\leq i\leq d} \\
& \xrightarrow[n \to +\infty]{a.s.} -E\left(\frac{\partial^2}{\partial\theta\partial\theta'}\ell_0(\theta_0)\right)
= E\left[\frac{1}{f_{\theta_0}^0}\left(\ddt f_{\theta_0}^0\right)\left(\ddt f_{\theta_0}^0\right)'\right] = \Sigma.
\end{align*}
\item[(c)]  If $U$ is a non-zero vector of $\R^d$, according to assumption Var, it holds that $U\ddt f_{\theta_0}^0 \neq 0$ a.s.
  Hence
  \[U\Sigma U' = E\left(\frac{1}{f_{\theta_0}^0}U\left(\ddt f^0_{\theta_0}\right)\left(\ddt f^0_{\theta_0}\right)'U'\right) \gr 0.\]
  Thus $\Sigma$ is positive definite.
\end{enumerate}

 \noindent From (a), apply the central limit theorem for stationary ergodic martingable difference
  sequence, it follows that
  \begin{equation}\label{tcl1_lt}
  \frac{1}{\sqrt{n}}\ddt L_n(T_{1,n},\theta_0) = \frac{1}{\sqrt{n}} \sum_{t=1}^n\ddt \ell_t (\theta_0) \xrightarrow[n\to+\infty]{\mathcal{D}} \mathcal{N} \left(0,E\left[\left(\ddt \ell_0(\theta_0)\right)\left(\ddt \ell_0(\theta_0)\right)'\right]\right).
  \end{equation}
  Recall that for $i = 1,\dots,d$, $\ddi \ell_t(\theta) = \left(\frac{Y_t}{f_{\theta}^t}-1\right)\ddi \ftt$.
 For $1 \leq i,j \leq d$, we have
 \begin{align*}
 E\Big(\ddi \ell_t(\theta_0)\times \ddj \ell_t(\theta_0)\Big) &=
 E\Big[E\Big(\big(\frac{Y_t}{f_{\theta_0}^t}-1\big)^2\ddi f_{\theta_0}^t \times \ddj f_{\theta_0}^t \big| \mathcal{F}_{t-1}\Big)\Big] \\
 &= E\Big[E\Big(\big(\frac{Y_t}{f_{\theta_0}^t}-1\big)^2 \big| \mathcal{F}_{t-1}\Big)\times \ddi f_{\theta_0}^t \times \ddj f_{\theta_0}^t \Big]
 \end{align*}
 But, we have,
 \begin{align*}
 E\Big[\Big(\frac{Y_t}{f_{\theta_0}^t}-1\Big)^2 \big| \mathcal{F}_{t-1}\Big]
 &= \frac{1}{(f_{\theta_0}^t)^2}E(Y_t^2 | \mathcal{F}_{t-1}) - \frac{2}{f_{\theta_0}^t}\times f_{\theta_0}^t + 1
 = \frac{1}{(f_{\theta_0}^t)^2}E(Y_t^2  | \mathcal{F}_{t-1}) - 1 \\
 &= \frac{1}{(f_{\theta_0}^t)^2} \Big(\var(Y_t | \mathcal{F}_{t-1}) + (E(Y_t | \mathcal{F}_{t-1}))^2 \Big) - 1
  = \frac{1}{(f_{\theta_0}^t)^2}( f_{\theta_0}^t + (f_{\theta_0}^t)^2) - 1 = \frac{1}{f_{\theta_0}^t}.
 \end{align*}
 It comes that,
 \[E\Big[\Big(\ddi \ell_t(\theta_0)\Big)\times \Big(\ddj \ell_t(\theta_0)\Big)'\Big]  =
 E\Big[\frac{1}{f_{\theta_0}^t}\Big(\ddi f_{\theta_0}^t\Big) \times \Big(\ddj f_{\theta_0}^t\Big)'\Big]. \]
 Hence
 \[E\Big[\Big(\ddt \ell_t(\theta_0)\Big)\times \Big(\ddt \ell_t(\theta_0)\Big)'\Big]  =
 E\Big[\frac{1}{f_{\theta_0}^t}\Big(\ddt f_{\theta_0}^t\Big) \times \Big(\ddt f_{\theta_0}^t\Big)'\Big] = \Sigma. \]
 Thus, (\ref{tcl1_lt}) becomes
\begin{equation}\label{tcl2_lt}
\frac{1}{\sqrt{n}}\ddt L_n(T_{1,n},\theta_0) = \frac{1}{\sqrt{n}}\sum_{t=1}^n \ddt \ell_t(\theta_0) \xrightarrow[n\to+\infty]{\mathcal{D}} \mathcal{N} (0,\Sigma).
\end{equation}
(b) and (c) implies that the matrix $G_n(T_{1,n},\tnc(T_{1,n}))$ converges a.s. to $\Sigma$ and $G_n(T_{1,n},\tnc(T_{1,n}))$ is invertible for $n$ large enough.
Hence, from (\ref{approx1_G}) and (\ref{tcl2_lt}), we have
\begin{align*}
\sqrt{n}(\tnc(T_{1,n})-\theta_0) &= \frac{1}{\sqrt{n}}\left(G_n(T_{1,n},\tnc(T_{1,n}))\right)^{-1}\ddt L_n(T_{1,n},\theta_0) + o_P(1) \\
 &= \frac{1}{\sqrt{n}}\Sigma^{-1}\ddt L_n(T_{1,n},\theta_0) + o_P(1)\\
& ~ \xrightarrow[n\to +\infty]{\mathcal{D}} \mathcal{N}(0,\Sigma^{-1}).
\end{align*}
                                                                                                                                   \Box
 ~ \\

   Before proving the Theorem \ref{theo3}, let us prove first some preliminary lemma.
   Under H$_0$, recall
\[\Sigma = E\Big[\frac{1}{f^0_{\theta_0}}\Big(\frac{\partial}{\partial\theta}f^0_{\theta_0}\Big)\Big(\frac{\partial}{\partial\theta}f^0_{\theta_0}\Big)'\Big]
= E\Big[\Big(\frac{\partial}{\partial\theta}\ell_0(\theta_0)\Big)\Big(\frac{\partial}{\partial\theta}\ell_0(\theta_0)\Big)'\Big].\]
Define the statistics
\begin{itemize}
  \item $\displaystyle C_n = \max_{v_n \leq k \leq n-v_n}C_{n,k}$ where
 \[C_{n,k} = \frac{1}{q^2\left(\frac{k}{n}\right)}\frac{k^2(n-k)^2}{n^3}\left(\widehat{\theta}_n(T_{1,k})-\widehat{\theta}_n(T_{k+1,n})\right)'
\Sigma\left(\widehat{\theta}_n(T_{1,k})-\widehat{\theta}_n(T_{k+1,n})\right);\]
  \item
$\displaystyle Q_n^{(1)} = \max_{v_n \leq k \leq
n-v_n}Q_{n,k}^{(1)}$ where
\[Q_{n,k}^{(1)} = \frac{k^2}{n}\left(\widehat{\theta}_n(T_{1,k})-\widehat{\theta}_n(T_{1,n})\right)'
\Sigma\left(\widehat{\theta}_n(T_{1,k})-\widehat{\theta}_n(T_{1,n})\right);\]
  \item
$\displaystyle Q_n^{(2)} = \max_{v_n \leq k \leq
n-v_n}Q_{n,k}^{(2)}$ where
\[Q_{n,k}^{(2)} = \frac{(n-k)^2}{n}\left(\widehat{\theta}_n(T_{k+1,n})-\widehat{\theta}_n(T_{1,n})\right)'
\Sigma\left(\widehat{\theta}_n(T_{k+1,n})-\widehat{\theta}_n(T_{1,n})\right).\]
\end{itemize}

 \begin{lem}\label{lem2}
Under assumptions of Theorem \ref{theo3}, as $n \to +\infty$,
\begin{itemize}
  \item[(i)] $\max_{v_n \leq k \leq n-v_n}\abs{\widehat{C}_{n,k}-C_{n,k}} = o_P(1)$;
  \item[(ii)] for $j=1,2, \max_{v_n \leq k \leq n-v_n}\abs{\widehat{Q}_{n,k}^{(j)}-Q_{n,k}^{(j)} } = o_P(1)$.
\end{itemize}
\end{lem}

 \begin{dem}
 (i) For any $v_n\leq k \leq n-v_n$, we have as $n \rightarrow \infty$

   \begin{align*}
 \abs{\widehat{C}_{n,k}-C_{n,k}} &=\frac{1}{q^2\left(\frac{k}{n}\right)}\frac{k^2(n-k)^2}{n^3} \abs{ \left(\widehat{\theta}_n(T_{1,k})-\widehat{\theta}_n(T_{k+1,n})\right)'
\Big(\widehat{\Sigma}_n(u_n)-\Sigma \Big)\left(\widehat{\theta}_n(T_{1,k})-\widehat{\theta}_n(T_{k+1,n})\right)} \\
&\leq \frac{1}{q^2\left(\frac{k}{n}\right)}\frac{k^2(n-k)^2}{n^3}\norme{\widehat{\Sigma}_n(u_n)-\Sigma} \norme{\widehat{\theta}_n(T_{1,k})-\widehat{\theta}_n(T_{k+1,n})}^2\\
&\leq C \frac{1}{ q^2\left(\frac{k}{n}\right) }  \frac{k(n-k)}{n^2}  \norme{\widehat{\Sigma}_n(u_n)-\Sigma}
   \Big(\norme{\sqrt{k}(\widehat{\theta}_n(T_{1,k})-\theta_0)}^2 + \norme{\sqrt{n-k}(\widehat{\theta}_n(T_{k+1,n}) - \theta_0)}^2 \Big) \\
& \leq C \frac{1}{ q^2\left(\frac{k}{n}\right) }  \frac{k(n-k)}{n^2} o(1) O_P(1).
\end{align*}
 Thus, as $n\rightarrow \infty$, it holds that
 \begin{align*}
 \max_{ v_n \leq k \leq n-v_n}\abs{\widehat{C}_{n,k}-C_{n,k}} &\leq o_P(1) \max_{ v_n \leq k \leq n-v_n} \frac{1}{ q^2\left(\frac{k}{n}\right) } \frac{k(n-k)}{n^2} \\
&\leq o_P(1) \max_{ \frac{v_n}{n} \leq \frac{k}{n} \leq 1- \frac{v_n}{n}} \frac{1}{ q^2\left(\frac{k}{n}\right) } \frac{k}{n} (1-\frac{k}{n} )\\
&\leq o_P(1) \sup_{ 0<\tau<1 } \Big(\frac{ \sqrt{\tau(1-\tau)} }{ q(\tau) } \Big)^2 = o_P(1).
\end{align*}
The last equality above holds because $\sup_{ 0<\tau<1 } \frac{ \sqrt{\tau(1-\tau)} }{ q(\tau)} <\infty$; it is a consequence of the properties of the function $q$
when $I_{0,1}(q,c)$ is finite for some $c>0$.\\
(ii) Goes the same lines as in (i).
 \end{dem}
  \paragraph{Proof of Theorem \ref{theo3}}
\begin{enumerate}
  \item According to Lemma \ref{lem2}, it suffices to show that
  \[C_n \xrightarrow[n \to +\infty]{\mathscr{D}} \sup_{0 \pet \tau \pet 1}\frac{\norme{W_d(\tau)}^2}{q^2(\tau)}.\]
  Let $v_n \leq k \leq n-v_n$. By applying (\ref{Taylor_G_2}) with $T = T_{k+1,n}$, we have
\[G_n(T_{1,k}, \tnc(T_{1,k}))(\tnc(T_{1,k})-\theta_0) = \frac{1}{k}\left( \ddt L_n(T_{1,k},\theta_0) - \ddt L_n(T_{1,k},\tnc(T_{1,k}))\right)\]
and
\begin{align*}
G_n(T_{k+1,n}, \tnc(T_{k+1,n}))(\tnc(T_{k+1,n})-\theta_0) = \frac{1}{n-k}\left(  \ddt L_n(T_{k+1,n},\theta_0)
 -\ddt L_n(T_{k+1,n},\tnc(T_{k+1,n}))   \right).
\end{align*}
As $n \to +\infty$, we have
\begin{align*}
& \norme{G_n(T_{1,k}, \tnc(T_{1,k})) - \Sigma} = \norme{G_n(T_{k+1,n},\tnc(T_{k+1,n}))-\Sigma} = o(1) \\
& \sqrt{k}(\tnc(T_{1,k})-\theta_0) = O_P(1) \esp{and} \sqrt{n-k}(\tnc(T_{k+1,n})-\theta_0) = O_P(1).
\end{align*}
Thus, we have
\begin{align*}
\sqrt{k}~\Sigma(\tnc(T_{1,k})-\theta_0) &= \frac{1}{\sqrt{k}}\left( \ddt L_n(T_{1,k},\theta_0) -\ddt L_n(T_{1,k},\tnc(T_{1,k}))  \right) \\
&\qquad -\sqrt{k}(G_n(T_{1,k},\tnc(T_{1,k}))-\Sigma)(\tnc(T_{1,k})-\theta_0) \\
&=\frac{1}{\sqrt{k}}\left(\ddt L_n(T_{1,k},\theta_0) - \ddt L_n(T_{1,k},\tnc(T_{1,k}))  \right) + o_P(1)
\end{align*}
i.~e.
\[
\Sigma(\tnc(T_{1,k})-\theta_0) =
\frac{1}{k}\left(\ddt L_n(T_{1,k},\theta_0) - \ddt L_n(T_{1,k},\tnc(T_{1,k})) \right) + o_P\left(\frac{1}{\sqrt{k}}\right).
\]
Moreover we have (as $n \to +\infty$)
\begin{align*}
\frac{1}{\sqrt{k}}\norme{\ddt L_n(T_{1,k},\tnc(T_{1,k}))} &= \frac{1}{\sqrt{k}}\norme{\ddt L_n(T_{1,k},\tnc(T_{1,k}))-0} \\
&\leq \frac{1}{\sqrt{k}}\norme{\ddt L_n(T_{1,k},\tnc(T_{1,k}))-\ddt \widehat{L}_n(T_{1,k},\tnc(T_{1,k}))} \\
&\leq \frac{1}{\sqrt{k}}\norme{\ddt L_n(T_{1,k},\theta)-\ddt \widehat{L}_n(T_{1,k},\theta)}_{\Theta} = o_P(1).
\end{align*}
Hence
\begin{align*}
\Sigma(\tnc(T_{1,k})-\theta_0) &= \frac{1}{\sqrt{k}}\left(\frac{1}{k} \ddt L_n(T_{1,k},\theta_0) -\frac{1}{\sqrt{k}}\ddt L_n(T_{1,k},\tnc(T_{1,k}))  \right)
+ o_P\left(\frac{1}{\sqrt{k}}\right) \\
&= \frac{1}{k}\ddt L_n(T_{1,k},\theta_0) + \frac{1}{\sqrt{k}}o_P(1)   + o_ p\left(\frac{1}{\sqrt{k}}\right) \\
&= \frac{1}{k}\ddt L_n(T_{1,k},\theta_0) + o_P(\frac{1}{\sqrt{k}}).
\end{align*}
Thus we have
\[\Sigma(\tnc(T_{1,k})-\theta_0) = \frac{1}{k}\ddt L_n(T_{1,k},\theta_0) + o_P\left(\frac{1}{\sqrt{k}}\right).\]
By going similarly lines, we obtain
\[\Sigma(\tnc(T_{k+1,n})-\theta_0) = \frac{1}{n-k}\ddt L_n(T_{k+1,n},\theta_0) + o_P\left(\frac{1}{\sqrt{n-k}}\right).\]
By subtracting the two above equalities, it follows that
\begin{align*}
\Sigma(\tnc(T_{1,k})-\tnc(T_{k+1,n})) &= \frac{1}{k}\ddt L_n(T_{1,k},\theta_0) - \frac{1}{n-k}\ddt L_n(T_{k+1,n},\theta_0)
+ o_P \Big(\frac{1}{\sqrt{k}}+\frac{1}{\sqrt{n-k}}\Big) \\
&= \frac{1}{k}\ddt L_n(T_{1,k},\theta_0) - \frac{1}{n-k} \Big( \ddt L_n(T_{1,n},\theta_0) - \ddt L_n(T_{1,k},\theta_0) \Big) \\
& \qquad + o_P \Big(\frac{1}{\sqrt{k}}+\frac{1}{\sqrt{n-k}} \Big) \\
&= \frac{n}{k(n-k)}\Big(\ddt L_n(T_{1,k},\theta_0) - \frac{n}{k}\ddt L_n(T_{1,n},\theta_0) \Big) + o_P \Big(\frac{1}{\sqrt{k}}+\frac{1}{\sqrt{n-k}}\Big)
\end{align*}
i.~e.
\begin{align*}
\frac{k(n-k)}{n^{\frac{3}{2}}}\Sigma(\tnc(T_{1,k})-\tnc(T_{k+1,n})) & = \frac{1}{\sqrt{n}}\left(\ddt L_n(T_{1,k},\theta_0) - \frac{k}{n}\ddt L_n(T_{1,n},\theta_0)\right) \\
 & \hspace{5cm} + o_P \Big(\frac{\sqrt{k(n-k)}}{n}\frac{\sqrt{k}+\sqrt{n-k}}{\sqrt{n}}\Big) \\
 &= \frac{1}{\sqrt{n}}\left(\ddt L_n(T_{1,k},\theta_0) - \frac{k}{n}\ddt L_n(T_{1,n},\theta_0)\right) + o_P(1).
\end{align*}
Thus
\begin{equation}
 \frac{k(n-k)}{n^{\frac{3}{2}}}\Sigma^{-1/2}\Sigma(\tnc(T_{1,k})-\tnc(T_{k+1,n})) = \frac{\Sigma^{-1/2}}{\sqrt{n}}\left(\ddt L_n(T_{1,k},\theta_0) - \frac{k}{n}\ddt L_n(T_{1,n},\theta_0)\right)
+ o_P(1). \label{approx1_theta}
\end{equation}
For $0 \pet \tau \pet 1$, we have
\[\frac{1}{\sqrt{n}} \ddt L_n(T_{1,[n\tau]},\theta_0) = \frac{1}{\sqrt{n}} \sum_{t=1}^{[n\tau]} \ddt \ell_t(\theta_0) .\]
We have shown (see the proof of Theorem \ref{theo2}) that $\left(\ddt \ell_t(\theta_0,\mathcal{F}_t)\right)_{t \in \Z}$ is a stationary
ergodic square integrable difference process with covariance matrix $\Sigma$. By the Central limit theorem for the
martingale difference sequence (see Billingstey, 1968), it holds that
\begin{align*}
& \frac{1}{\sqrt{n}}\left(\ddt L_n(T_{1,[n\tau]}, \theta_0) - \frac{[n\tau]}{n}\ddt L_n(T_{1,n},\theta)\right) \\
& =\frac{1}{\sqrt{n}}\left(\sum_{t=1}^{[n\tau]}\ddt \ell_t(\theta_0)-\frac{[n\tau]}{n}\sum_{t=1}^n\ddt \ell_t(\theta_0)\right)
\xrightarrow[n\to +\infty]{ \mathcal{D} } B_{\Sigma}(\tau)-\tau B_{\Sigma}(1)
\end{align*}
where $B_{\Sigma}$ is a  Gaussian process with covariance matrix $\min(s,\tau)\Sigma$. Thus it follows that
\[
\frac{1}{\sqrt{n}}\Sigma^{-1/2}\left(\ddt L_n(T_{1,[n\tau]},\theta_0) - \frac{[n\tau]}{n}\ddt L_n(T_{1,n},\theta_0)\right)
\xrightarrow[n\to +\infty]{\mathcal{D}} B_d(\tau)-\tau B_d(1) = W_d(\tau)
\]
in $D([0,1])$, where $B_d$ is a $d$-dimensional standard motion, and $W_d$ is a $d$-dimensional
Brownian bridge.

So, we have (see (\ref{approx1_theta}))
\begin{align*}
C_{n,[n\tau]} &= \frac{[n\tau]^2(n-[n\tau])^2}{n^3}(\tnc(T_{1,[n\tau]})-\tnc(T_{[n\tau]+1,n}))'\Sigma(\tnc(T_{1,[n\tau]})-\tnc(T_{[n\tau]+1,n})) \\
& = \Big \| \frac{[n\tau](n-[n\tau])}{n^{3/2}}\Sigma^{-1/2}\Sigma(\tnc(T_{1,[n\tau]})-\tnc(T_{[n\tau]+1,n})) \Big \|^2 \\
& = \Big \| \frac{1}{\sqrt{n}} \Sigma^{-1/2} \Big( \sum_{t=1}^{[n\tau]}\ddt \ell_t(\theta_0)- \frac{[n\tau]}{n}\sum_{t=1}^n\ddt \ell_t(\theta_0) \Big) \Big \|^2 +
 o_P(1)\xrightarrow[n\to +\infty]{\mathcal{D}} \norme{W_d(\tau)}^2 \esp{in} \mathcal{D}([0,1]).
\end{align*}
 Hence, according to the properties  of $q$, we have for any $0 \pet \epsilon \pet 1/2$
\begin{align*}
\max_{[n\epsilon] \leq k \leq n-[n\epsilon]}C_{n,k} &=
\max_{[n\epsilon] \leq k \leq
n-[n\epsilon]}\frac{1}{q^2\left(\frac{k}{n}\right)}\frac{k^2(n-k)^2}{n^3}\left(\tnc(T_{1,k})-\tnc(T_{k+1,n})\right)'
\Sigma\left(\tnc(T_{1,k})-\tnc(T_{k+1,n})\right) \\
&= \sup_{\epsilon \leq \tau \leq 1-\epsilon}\frac{1}{q^2\left(\frac{[n\tau]}{n}\right)}\frac{[n\tau]^2(n-[n\tau])^2}{n^3}\left(\tnc(T_{1,[n\tau]})-\tnc(T_{[n\tau]+1,n})\right)'\times \\
& \qquad \Sigma\left(\tnc(T_{1,[n\tau]})-\tnc(T_{[n\tau]+1,n})\right) \\
&= \sup_{\epsilon \leq \tau \leq 1-\epsilon} \Big \| \frac{1}{q\left(\frac{[n\tau]}{n}\right)}\frac{[n\tau](n-[n\tau])}{n^{3/2}}
\Sigma^{1/2}\left(\tnc(T_{1,[n\tau]})-\tnc(T_{[n\tau]+1,n})\right)\Big \|^2 \\
&= \sup_{\epsilon \leq \tau \leq 1-\epsilon} \Big \|\frac{1}{q\left(\frac{[n\tau]}{n}\right)}\frac{1}{\sqrt{n}}
\left(\sum_{t=1}^{[n\tau]}\ddt \ell_t(\theta_0)-\frac{[n\tau]}{n}\sum_{t=1}^n\ddt \ell_t(\theta_0)\right) \Big \|^2 + o_P(1) \\
& \xrightarrow[n \to+\infty]{\mathcal{D}} \sup_{\epsilon \leq \tau \leq 1-\epsilon}\frac{\norme{W_d(\tau)}^2}{q^2(\tau)}.
\end{align*}
Therefore, we have shown that
\[C_{n,[n\tau]}\xrightarrow[n\to+\infty]{\mathcal{D}} \frac{\norme{W_d(\tau)}^2}{q^2(\tau)} \esp{in} \mathcal{D}([0,1])\]
and that for all $\epsilon \in (0,1/2)$,
\[\max_{[n\epsilon]\leq k \leq n-[n\epsilon]} C_{n,k} = \sup_{\epsilon \leq \tau \leq 1-\epsilon}C_{n,[n\tau]}
\xrightarrow[n \to+\infty]{\mathcal{D}} \sup_{\epsilon \leq \tau \leq 1-\epsilon}\frac{\norme{W_d(\tau)}^2}{q^2(\tau)}.\]
%
%
Moreover, since $I_{0,1}(q,c) \pet +\infty$ for some $c \gr 0$, one can show that almost surely $\lim_{\tau \to 0}\frac{\norme{W_d(\tau)}}{q(\tau)} \pet \infty$
and $\lim_{\tau \to 1}\frac{\norme{W_d(\tau)}}{q(\tau)} \pet \infty$ (see for instance \cite{Csorgo1986}). Hence, for $n$ large enough we have
\[C_n = \max_{v_n \leq k \leq n-v_n} C_{n,k} = \sup_{\frac{v_n}{n} \leq \tau \leq 1-\frac{v_n}{n}} C_{n,[n\tau]}
\xrightarrow[n \to+\infty]{\mathcal{D}} \sup_{0 \leq \tau \leq 1}\frac{\norme{W_d(\tau)}^2}{q^2(\tau)}.\]
\item Apply Lemma \ref{lem2} and goes along similar lines as in the proof of Theorem 1 of \cite{Kengne2011}.  \Box
\end{enumerate}
 \paragraph{Proof of Theorem~\ref{theo4}}~

\begin{enumerate}
  \item Assume the alternative with one change at $t_1^* = [\tau_1^*n]$ with $0 \pet \tau_1^* \pet 1$.
  The observation satisfy
  \[Y_t =
  \begin{cases}
    Y_t^{(1)} & \text{for} \quad t \leq t_1^*, \\
    Y_t^{(2)} & \text{for} \quad t \gr t_1^*.
  \end{cases}
  \]
 Where $(Y_t^{(1)})$ and $Y_t^{(2)}$ satisfy the main model~$(3)$
 i.~e.
 \[\left.Y_t^{(i)}\right/\mathcal{F}_{t-1}^{(i)} \sim \mathcal{P}(\lambda_t^{(i)}) \esp{with} \lambda_t^{(i)} = f_{\theta_i^*}(Y_{t-1}^{(i)}, Y_{t-2}^{(i)}, \dots ), \theta_1^*
 \neq \theta_2^*, \esp{and} \mathcal{F}_t^{(i)} = \sigma(Y_s^{(i)}, s\leq t).\]

 Recall that $\widehat{C}_n = \max_{v_n \leq k \leq n-v_n}\widehat{C}_{n,k} \geq \widehat{C}_{n,t_1^*}$. It suffices to
 show that $\widehat{C}_{n,t_1^*} \xrightarrow[n\to +\infty]{P}+\infty$.

 We have
 \[\widehat{C}_{n,t_1^*} = \frac{1}{q\left(\frac{t_1^*}{n}\right)}\frac{{t_1^*}^2(n-t_1^*)^2}{n^3}\left(\tnc(T_{1,t_1^*})-\tnc(T_{t_1^*+1,n})\right)'
 \widehat{\Sigma}(u_n)\left(\tnc(T_{1,t_1^*})-\tnc(T_{t_1^*+1,n})\right).\]
 Recall that the likelihood function computed on any subset $T \subset \{1, \dots, n\}$ is defined by
 \[\widehat{L}_n(T,\theta) = \sum_{t\in T}\widehat{\ell}_t(\theta),\esp{where} \widehat{\ell}_t(\theta) = Y_t\log\widehat{\ell}_{\theta}^t - \widehat{\ell}_{\theta}^t, \esp{with}
 \widehat{f}_{\theta}^t = f_{\theta}(Y_{t-1}, Y_{t-2},\dots, Y_1, 0, \dots).\]

 So, for any $t \in \{1, \dots, t_1^*\}$, $\widehat{f}_{\theta}^t = f_{\theta}(Y_{t-1}^{(1)}, Y_{t-2}^{(1)},\dots, Y_1^{(1)}, 0,
 \dots)$. Then $\theta \mapsto \widehat{L}(T_{1,t_1^*},\theta)$ is the likelihood function of the stationary process $(Y_t^{(1)})_{t \in \Z}$
 computed on $\{1, \dots, t_1^*\}$. According to Theorem \ref{theo1}, it holds that $\tnc(T_{1,t_1^*}) \xrightarrow[n \to +\infty]{a.s.}\theta_1^*$.
 Moreover, recall that
 \begin{align*}
 \widehat{\Sigma}_n(u_n) &= \frac{1}{2}\left[\left.\left(\frac{1}{u_n}\sum_{t=1}^{u_n}\frac{1}{\fttc}\left(\ddt\fttc\right)\left(\ddt\fttc\right)'\right)\right|_{\theta = \tnc(T_{1,u_n})} \right. \\
 &\qquad \qquad +\left.\left.\left(\frac{1}{n-u_n}\sum_{t=u_n+1}^{n}\frac{1}{\fttc}\left(\ddt\fttc\right)\left(\ddt\fttc\right)'\right)\right|_{\theta = \tnc(T_{u_n+1,n})} \right].
 \end{align*}

 Denote $l_{t,1}(\theta) = Y_t^{(1)}\log f_{\theta}^{t,1} - f_{\theta}^{t,1}$ where
 \[f_{\theta}^{t,1} = f_{\theta}(Y_{t-1}^{(1)}, Y_{t-2}^{(1)}, \dots, Y_{1}^{(1)}, 0\dots).\] We have
 $\tnc (T_{1,u_n}) \xrightarrow[n \to +\infty]{a.s.}\theta_1^*$ and, from Lemma \ref{lem1}, it holds that
 \[\frac{1}{u_n}\sum_{t=1}^{u_n}\frac{1}{\fttc}\left(\ddt\fttc\right)\left(\ddt\fttc\right)' \xrightarrow[n\to +\infty]{a.s.}\Sigma^{(1)} \esp{where}
 \Sigma^{(1)} = E\left[\frac{1}{f_{\theta_1^*}^{0,1}}\left(\ddt f_{\theta_1^*}^{0,1}\right)\left(\ddt f_{\theta_1^*}^{0,1}\right)'\right]. \]
 We have
 \[\widehat{L}_n(T_{t_1^*+1,n},\theta) = \sum_{t=t_1^*+1}^n\widehat{\ell}_t(\theta) \esp{where} \widehat{\ell}_t(\theta) = Y_t\log\fttc - \fttc, \]
 with
 \[\fttc = f_{\theta}(Y_{t-1}, Y_{t-2}, \dots, Y_1, 0, \dots) = f_{\theta}(Y_{t-1}^{(1)}, Y_{t-2}^{(1)}, \dots, Y_1^{(1)}, 0, \dots).\]
 Define
 \[\widehat{L}_{n,2}(T_{t_1^*+1,n},\theta) = \sum_{t=t_1^*+1}^n\widehat{\ell}_{t,2}(\theta), \esp{where} \widehat{\ell}_{t,2}(\theta) = Y_t^{(2)}\log \widehat{f}_{\theta}^{t,2} - \widehat{f}_{\theta}^{t,2},\]
 with $\widehat{f}_{\theta}^{t,2} = f_{\theta}(Y_{t-1}^{(2)}, Y_{t-2}^{(2)}, \dots, Y_{1}^{(2)},0,\dots)$.

 Remarks that the difference between $\fttc$ and $\widehat{f}_{\theta}^{t,2}$ lies on the dependence with the past.
 $\fttc$ can contain $Y_{t-1}^{(1)}$, but not $\widehat{f}_{\theta}^{t,2}$. $\theta \mapsto
\widehat{L}_{n,2}(T_{t_1^*+1,n},\theta)$ is the approximated likelihood
of the stationary model after change. By Theorem \ref{theo1}, it
holds that
\[\widehat{\theta}_n^{(2)}(T_{t_1^*+1,n}) = \argmax_{\theta\in \Theta} \widehat{L}_{n,2}(T_{t_1^*+1,n},\theta) \xrightarrow[+ \to +\infty]{a.s.}\theta_2^*.\]
Let us show that
\[\frac{1}{n-t_1^*}\norme{\widehat{L}_n(T_{t_1^*+1,n},\theta)-\widehat{L}_{n,2}(T_{t_1^*+1,n},\theta)}_{\Theta}\xrightarrow[n \to +\infty]{a.s.} 0.\]
We have, as
\[\frac{1}{n-t_1^*}\norme{\widehat{L}_n(T_{t_1^*+1,n},\theta)-\widehat{L}_{n,2}(T_{t_1^*+1,n},\theta)}_{\Theta} \leq
\frac{1}{n-t_1^*}\sum_{k=1}^{n-t_1^*}\norme{\widehat{\ell}_{t_1^*+k}(\theta)-\widehat{\ell}_{t_1^*+k,2}(\theta)}_{\Theta},\]
using again Kounias ($1969$), it suffices to show that
\[\sum_{k \geq 1}\frac{1}{k}\E \norme{\widehat{\ell}_{t_1^*+k}(\theta), \widehat{\ell}_{t_1^*+k,2}(\theta)}_{\Theta}
= \sum_{k \geq 1} \frac{1}{t-t_1^*} \E \norme{\widehat{\ell}_t(\theta)-\widehat{\ell}_{t,2}(\theta)}_{\Theta} \pet \infty. \]
 For $t \geq t_1^*+1$, we have
 \begin{align*}
 E\left[\norme{\widehat{\ell}_t(\theta)-\widehat{\ell}_{t,2}(\theta)}_{\Theta}\right] &= E\left[\norme{Y_t^{(2)}\log \fttc-\fttc-Y_t^{(2)}\log \widehat{f}_{\theta}^{t,2}+\widehat{f}_{\theta}^{t,2}}_{\Theta}\right] \\
 &\leq E\left[\abs{Y_t^{(2)}}\norme{\log \fttc-\log
 \widehat{f}_{\theta}^{t,2}}_{\Theta} + \norme{\fttc -
 \widehat{f}_{\theta}^{t,2}}_{\Theta}\right].
 \end{align*}
 We can show, as in the proof of Theorem \ref{theo1}, that $\norme{\log \fttc-\log \widehat{f}_{\theta}^{t,2}}_{\Theta} \leq \frac{1}{\nc}\norme{\fttc - \widehat{f}_{\theta}^{t,2}}_{\Theta}$. Hence
 \begin{align*}
 E\Big[\norme{\widehat{\ell}_t(\theta)-\widehat{\ell}_{t,2}(\theta)}_{\Theta}\Big] &\leq E\Big[\Big(\frac{|Y_t^{(2)}|}{\nc}+1\Big)
 \norme{\fttc-\widehat{f}_{\theta}^{t,2}}_{\Theta}\Big] \\
 &\leq \Big(\E\Big[\Big(\frac{|Y_t^{(2)}|}{\nc}+1\Big)^2\Big]\Big)^{1/2}\times \left(E \norme{ \fttc - \widehat{f}_{\theta}^{t,2} }_{\Theta}^2  \right)^{1/2} \\
 &\leq C\left(E \norme{ \fttc - \widehat{f}_{\theta}^{t,2} }_{\Theta}^2  \right)^{1/2}.
 \end{align*}
 But, for $t \geq t_1^*+1$, we have
 \begin{align*}
 \norme{\fttc-\widehat{f}_{\theta}^{t,2}}_{\Theta} &=
 \norme{f_{\theta}(Y_{t-1}^{(2)}, \dots, Y_{t_1^*+1}^{(2)},
 Y_{t_1^*}^{(1)}, \dots,Y_{1}^{(1)},0,\dots) -
 f_{\theta}(Y_{t-1}^{(2)}, \dots, Y_{t_1^*+1}^{(2)}, 0,\dots)} \\
 &\leq \sum_{j=t-t_1^*+1}^t \alpha_j^{(0)}(\Theta)|Y_{t-j}^{(1)}|.
\end{align*}
Thus, by using Minkowski inequality, it holds that
\begin{align*}
E\big[\norme{\widehat{\ell}_t(\theta)-\widehat{\ell}_{t,2}(\theta)}_{\Theta}\big] &\leq C \big(E\norme{\fttc - \widehat{f}_{\theta}^{t,2}}_{\Theta}^2\big)^{1/2} \\
&\leq C\Big(E\Big[\Big(\sum_{j=t-t_1^*+1}^t \alpha_j^{(0)}(\Theta)|Y_{t-j}^{(1)}|\Big)^2\Big]\Big)^{1/2} \\
&\leq \sum_{j=t-t_1^*+1}^t \alpha_j^{(0)}(\Theta)\Big(E\Big[|Y_{t-j}^{(1)}|^2\Big]\Big)^{1/2} \\
&\leq C\sum_{j=t-t_1^*+1}^t \alpha_j^{(0)}(\Theta).
\end{align*}
Thus, it comes that
\begin{align*}
\sum_{t\geq t_1^*+1}\frac{1}{t-t_1^*}\norme{\widehat{\ell}_t(\theta)-\widehat{\ell}_{t,2}(\theta)}_{\Theta} &\leq
C\sum_{t\geq t_1^*+1}\frac{1}{t-t_1^*}\sum_{j=t-t_1^*+1}^t \alpha_j^{(0)}(\Theta) \\
&\leq \sum_{t\geq t_1^*+1}\sum_{j=t-t_1^*+1}^t\frac{1}{t-t_1^*}\alpha_j^{(0)}(\Theta).
\end{align*}
Set $l = t-t_1^*$, we have
\begin{align*}
\sum_{t\geq t_1^*+1}\sum_{j=t-t_1^*+1}^t\frac{1}{t-t_1^*}\alpha_j^{(0)}(\Theta) &=
\sum_{l \geq 1}\sum_{j=l+1}^{l+t_1^*}\frac{1}{l}\alpha_j^{(0)}(\Theta) \\
&\leq \sum_{l\geq 1}\sum_{j=l}^{l+t_1^*}\frac{1}{l}\alpha_j^{(0)}(\Theta) \\
&= \sum_{j=1}^{t_1^*}\sum_{l=1}^j\frac{1}{l}\alpha_j^{(0)}(\Theta) + \sum_{j=t_1^*+1}^{+\infty}\sum_{l=j-t_1^*}^j\frac{1}{l}\alpha_j^{(0)}(\Theta) \\
&\leq C\sum_{j=1}^{t_1^*}\alpha_j^{(0)}(\Theta)\log j + \sum_{j=t_1^*+1}^{+\infty}\alpha_j^{(0)}(\Theta)\sum_{l=1}^j\frac{1}{l} \\
&\leq C\sum_{j\geq 1}\alpha_j^{(0)}(\Theta)\log j \pet +\infty.
\end{align*}
So, we have
\[\sum_{t\geq t_1^*+1}\frac{1}{t-t_1^*}\norme{\widehat{\ell}_t(\theta)-\widehat{\ell}_{t,2}(\theta)}_{\Theta} \leq C \pet \infty\esp{a.s.}\]
Hence,
\[\frac{1}{n-t_1^*}\norme{\widehat{L}_n(T_{t_1^*,n},\theta)-\widehat{L}_{n,2}(T_{t_1^*,n},\theta)}_{\Theta}
\xrightarrow[n \to+\infty]{a.s.} 0.\]
According to the proof of Theorem \ref{theo1},
\[\frac{1}{n-t_1^*}\norme{\widehat{L}_{n,2}(T_{t_1^*,n},\theta)-\E \ell_{0,2}(\theta)}_{\Theta} \xrightarrow[n \to+\infty]{a.s.} 0\]
%
%
where
 \[l_{t,2}(\theta) = Y_t^{(2)}\log f_{\theta}^{t,2}-f_{\theta}^{t,2}, \esp{with} f_{\theta}^{t,2} = f_{\theta}(Y_{t-1}^{(2)}, Y_{t-2}^{(2)}, \dots).\]
 Moreover, the function $\theta \mapsto \E (l_{0,2}(\theta))$ has a unique maximum at $\theta_2^*$. This is enough to conclude that
 $\tnc(T_{t_1^*+1,n})\xrightarrow[n \to +\infty]{a.s.} \theta_2^*$.
 To complete the proof of this part of Theorem \ref{theo4}, remarks that the two matrices in the definition of $\widehat{\Sigma}_n(u_n)$ are
 positive semi-definite (by definition) and the first one converges a.s. to $\Sigma^{(1)}$ which is positive definite. Thus, for large enough, we have
 \begin{align*}
 \widehat{C}_n &= \max_{v_n \leq k \leq n-v_n}\widehat{C}_{n,k} \geq \widehat{C}_{n,t_1^*}\\
 &\geq \frac{1}{q\left(\frac{t_1^*}{n}\right)}\frac{t_1^{* 2}(n-t_1^*)^2}{n^3}\left(\tnc(T_{1,t_1^*})-\tnc(T_{t_1^*+1},n)\right)' \times \\
 &  \qquad \left[\frac{1}{2}\left.\left(\frac{1}{u_n}\sum_{t=1}^{u_n}\frac{1}{\fttc}\left(\ddt \fttc\right)\left(\ddt \fttc\right)'\right)\right|_{\theta=\tnc(T_{1,u_n})}\right]
 \times \left(\tnc(T_{1,t_1^*})-\tnc(T_{t_1^*+1,n})\right)
 \end{align*}
\begin{align*}
 &\geq \frac{1}{\sup_{0\pet \tau \pet 1}q(\tau)}n(\tau_1^*(1-\tau_1^*))^2\left(\tnc(T_{1,t_1^*})-\tnc(T_{t_1^*+1,n})\right)'\times \\
 &\qquad \left[\frac{1}{2}\left.\left(\frac{1}{u_n}\sum_{t=1}^{u_n}\frac{1}{\fttc}\left(\ddt \fttc\right)\left(\ddt \fttc\right)'\right)\right|_{\theta=\tnc(T_{1,u_n})}\right]
 \times \left(\tnc(T_{1,t_1^*})-\tnc(T_{t_1^*+1,n})\right) \\
 &\geq C \times n\left(\tnc(T_{1,t_1^*})-\tnc(T_{t_1^*+1,n})\right)
 \times \left[\frac{1}{2}\left.\left(\frac{1}{u_n}\sum_{t=1}^{u_n}\frac{1}{\fttc}\left(\ddt \fttc\right)\left(\ddt \fttc\right)'\right)\right|_{\theta=\tnc(T_{1,u_n})}\right] \times \\
 & \qquad \left(\tnc(T_{1,t_1^*})-\tnc(T_{t_1^*+1,n})\right)
 \xrightarrow[n \to+\infty]{a.s.} + \infty.
 \end{align*}
 This holds because $\tnc(T_{1,t_1^*})-\tnc(T_{t_1^*+1,n}) \xrightarrow[n \to+\infty]{a.s.} \theta_1^{*}-\theta_2^* \neq 0$, and
 \[\frac{1}{u_n}\sum_{t=1}^{u_n}\frac{1}{\fttc}\left(\ddt \fttc\right)\left(\ddt \fttc\right)'\xrightarrow[n \to+\infty]{a.s.}\Sigma^{(1)},\]
 which is positive definite.

 This completes the first part of the proof of Theorem~\ref{theo4}.

 \item It goes along the same line as in the proof of Theorem 2 of \cite{Kengne2011}, by using the approximation of likelihood
 as above. \Box
 \end{enumerate}

 \section*{ Acknowledgements}
  The authors thank Konstantinos Fokianos for his support, especially for his help  in the real data application study.

\end{document}